\documentstyle[12pt,toc_small,bib_toc,twoside]{article}

\input{psfig}
 \def\tty{\tt}

\newfont{\bbf}{msbm10 at 12pt}
\def\ttt{\tty}

\def\zz{\mbox{\bbf Z}}
\def\rz{\mbox{\bbf R}}
\def\cz{\mbox{\bbf C}}

\def\disk{\mbox{\,\bbf D}}

\def\Circle{\mbox{\bbf S}^1}

\def\K{ {\mbox{\ttt K}} }

\def\*{ {\mbox{\tt $\star$}} }
\def\0{ {\mbox{\tt 0}} }
\def\1{ {\mbox{\tt 1}} }
\def\2{ {\mbox{\tt 2}} }
\def\3{ {\mbox{\tt 3}} }
\def\4{ {\mbox{\tt 4}} }

\def\ovl{\overline}

\def\<{\prec}
\def\>{\succ}

\def\phi1{\phi}
\def\phi{\varphi}

\ifx\ThmNum\undefined 
   \newtheorem{theorem}{Theorem}[section]
   \def\newsection#1{\setcounter{theorem}{0} \section{#1}}
\else 
   \newtheorem{theorem}{Theorem} 
   \def\newsection#1{\section{#1}}
\fi

\newtheorem{proposition}[theorem]{Proposition}
\newtheorem{lemma}[theorem]{Lemma}
\newtheorem{definition}[theorem]{Definition}

\newtheorem{corollary}[theorem]{Corollary}

\def\proof{\par\medskip\noindent {\sc Proof. }}
\def\proofof #1 {\par\medskip\noindent {\sc Proof of #1. }}

\def\sketchof #1 {\par\medskip\noindent {\sc Sketch of proof of #1. }}
\def\Box{\framebox[10pt]{\rule{0pt}{3pt}}}
\def\qed{\hfill $\Box$ \medskip \par}
\def\qedd{\hfill $\Box$}
\def\remark{\par\medskip \noindent {\sc Remark. }}
\def\lineclear{\rule{0pt}{0pt}\par\noindent}

\newfont{\sfH}{cmss12}
\def\sf{\sfH}
\def\reminder #1 {{\sf #1}}
\def\hide #1 {}

\def\SBIMSMark#1#2#3{
 \font\SBF=cmss10 at 10 true pt
 \font\SBI=cmssi10 at 10 true pt
 \setbox0=\hbox{\SBF Stony Brook IMS Preprint \##1}
 \setbox2=\hbox to \wd0{\hfil \SBI #2}
 \setbox4=\hbox to \wd0{\hfil \SBI #3}
 \setbox6=\hbox to \wd0{\hss
             \vbox{\hsize=\wd0 \parskip=0pt \baselineskip=10 true pt
                   \copy0 \break%
                   \copy2 \break%
                   \copy4 \break}}
 \dimen0=\ht6   \advance\dimen0 by \vsize \advance\dimen0 by 8 true pt
                \advance\dimen0 by -\pagetotal
 \dimen2=\hsize \advance\dimen2 by .25 true in
  \openin2=publishd.tex
  \ifeof2\setbox0=\hbox to 0pt{}
  \else 
     \setbox0=\hbox to 3.1 true in{
                \vbox to \ht6{\hsize=3 true in \parskip=0pt  \noindent  
                {\SBI Published in modified form:}\hfil\break
                \input publishd.tex 
                \vfill}}
  \fi
  \closein2
  \ht0=0pt \dp0=0pt
 \ht6=0pt \dp6=0pt
 \setbox8=\vbox to \dimen0{\vfill \hbox to \dimen2{\copy0 \hss \copy6}}
 \ht8=0pt \dp8=0pt \wd8=0pt
 \copy8
 \message{*** Stony Brook IMS Preprint #1, #2 ***}
}

\def\M{{\bf M}}
\def\Part#1{{\cal P}_{#1}}
\edef\theta{\vartheta}
\newcommand{\brkOK}{\discretionary{}{}{}}

\newlength\captionwidth
\def\newsection#1{
   \setcounter{theorem}{0} \section{#1}
   \markboth{#1}{Rational Parameter Rays of the Mandelbrot Set}
}

\title{\vspace{-15mm}Rational Parameter Rays\\ of the Mandelbrot Set}
\author{Dierk Schleicher}

\begin{document}
\thispagestyle{empty}

\maketitle
\SBIMSMark{1997/13}{November 1997}{}
\thispagestyle{empty}


\begin{center}
\begin{minipage}{80mm}
\def\toc_vspace{1.5em}
\def\tocname{Contents}
\tableofcontents{0.3em}
\end{minipage}
\end{center}
\vspace{5mm}

\begin{abstract}
We give a new proof that all external rays of the Mandelbrot set at
rational angles land, and of the relation between the external angle
of such a ray and the dynamics at the landing point. Our proof is
different from the original one, given by Douady and Hubbard and
refined by Lavaurs, in several ways: it replaces analytic arguments by
combinatorial ones; it does not use complex analytic dependence of the
polynomials with respect to parameters and can thus be made to apply
for non-complex analytic parameter spaces; this proof is also
technically simpler. Finally, we derive several corollaries about
hyperbolic components of the Mandelbrot set.

Along the way, we introduce partitions of dynamical and parameter
planes which are of independent interest, and we interpret the
Mandelbrot set as a symbolic parameter space of kneading sequences
and internal addresses.
\bigskip

Nous donnons une nouvelle d\'emonstration que tous les rayons
externes \`a arguments rationels de l'ensemble Mandelbrot 
aboutissent, et nous montrons la relation entre l'argument externe
d'un tel rayon et la dynamique au param\`etre o\`u le rayon aboutit.
Notre d\'emonstration est diff\'erente de l'originale, donn\'ee par
Douady et Hubbard et elabor\'ee par Lavaurs, \`a plusieurs
\'egards: elle remplace des arguments analytiques par des arguments
combinatoires; elle n'utilise pas la d\'ependence analytique des
polyn\^omes par rapport au param\`etre et peut donc \^etre
appliqu\'ee aux espaces de param\`etres qui ne sont pas 
analytiques complexes; la d\'emonstration est aussi techniquement
plus facile. Finalement, nous d\'emontrons quelques corollaires sur
les composantes hyperboliques de l'ensemble Mandelbrot. 

En route, nous introduisons des partitions du plan dynamique et de
l'espace des param\`etres qui sont int\'eressantes en elles-m\^emes,
et nous interpr\'etons l'en\-sem\-ble Mandelbrot comme un espace de
param\`etres symboliques contenant des {\sl kneading sequences} et
des adresses internes.

\end{abstract}
\pagebreak

\newsection {Introduction}
\label{SecIntro}

Quadratic polynomials, when iterated, exhibit amazingly rich
dynamics. Up to affine conjugation, these polynomials can be
parametrized uniquely by a single complex variable. The Mandelbrot
set serves to organize the space of (conjugacy classes of) quadratic
polynomials. It can be understood as a ``table of contents'' to the
dynamical possibilities and has a most beautiful structure. Much of
this structure has been discovered and explained in the groundbreaking
work of Douady and Hubbard \cite{Orsay}, and a deeper understanding
of the fine structure of the Mandelbrot set is a very active area of
research. The importance of the Mandelbrot set is due to the fact
that it is the simplest non-trivial parameter space of analytic
families of iterated holomorphic maps, and because of its
universality as explained by Douady and Hubbard \cite{Polylike}: the
typical local configuration in one-dimensional complex parameter
spaces is the Mandelbrot set (see also \cite{CurtUniv}).

\begin{figure}[htbp]
\centerline{\psfig{figure=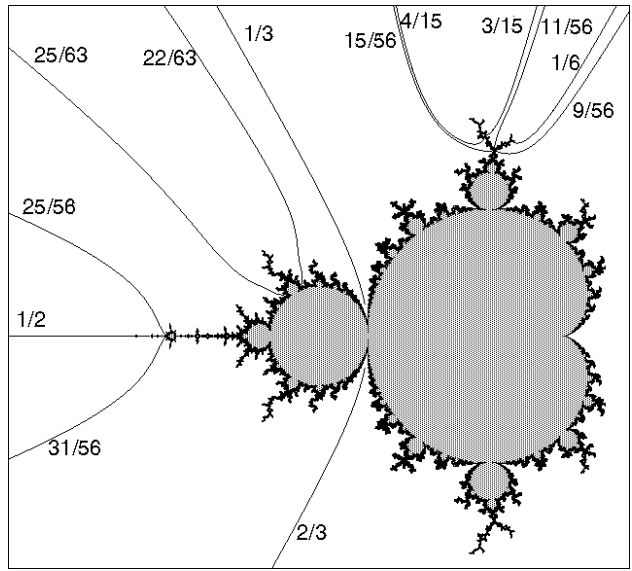,height=75mm}}
\centerline{\parbox{\captionwidth}{ 
\caption{\sl The Mandelbrot set and several of its parameter rays
which are mentioned in the text. (Picture courtesy of Jack Milnor)}
\label{FigMandelRays}
} }
\end{figure}

Unfortunately, most of the beautiful results of Douady and Hubbard on
the structure of the Mandelbrot set are written only in preliminary
form in the preprints \cite{Orsay}. The purpose of this article is
to provide concise proofs of several of their theorems. Our proofs
are quite different from the original ones in several respects: while
Douady and Hubbard used elaborate perturbation arguments for many
basic results, we introduce partitions of dynamical and parameter
planes, describe them by symbolic dynamics, and reduce many of the
questions to a combinatorial level. We feel that our proofs are
technically significantly simpler than those of Douady and Hubbard.
An important difference for certain applications is that our proof
does not use complex analytic dependence of the maps with respect to
the parameter and is therefore applicable in certain wider
circumstances: the initial motivation for this research was a project
with Nakane (see \cite{Nakane} and the references therein) to
understand the parameter space of antiholomorphic quadratic
polynomials, which depends only real-analytically on the parameter. Of
course, the ``standard proof'' using Fatou coordinates and Ecalle
cylinders, as developed by Douady and Hubbard and elaborated by
Lavaurs~\cite{LaEc}, is a most powerful tool giving interesting
insights; it has had many important applications. Our goal is to
present an alternative approach in order to enlarge the toolbox for
applications in different situations. 

The fundamental result we want to describe in this article is the
following theorem about landing properties of external rays of the
Mandelbrot set, a theorem due to Douady and Hubbard; for background
and terminology, see the next section.

\begin{theorem}[The Structure Theorem of the Mandelbrot Set]
\label{ThmRatRays} \lineclear
\vspace{ -24pt}
\begin{enumerate}
\item  
Every parameter ray at a periodic angle $\theta$ lands at a parabolic
parameter such that, in its dynamic plane, the dynamic ray at angle
$\theta$ lands at the parabolic orbit and is one of its two
characteristic rays.
\item  
Every parabolic parameter $c$ is the landing point of exactly
two parameter rays at periodic angles. These angles are the
characteristic angles of the parabolic orbit in the dynamic plane of
$c$. 
\item  
Every parameter ray at a preperiodic angle $\theta$ lands at a
Misiurewicz point such that, in its dynamic plane, the dynamic ray at
angle $\theta$ lands at the critical value.
\item  
Every Misiurewicz point $c$ is the landing point of a finite non-zero 
number of parameter rays at preperiodic angles. These angles are
exactly the external angles of the dynamic rays which land at the
critical value in the dynamic plane of $c$. 
\end{enumerate}
\end{theorem}

(The parameter $c=1/4$ is the landing point of a single parameter
ray, but this ray corresponds to external angles $0$ and $1$; we
count this ray twice in order to avoid having to state exceptions.)

The organization of this article is as follows: in
Section~\ref{SecDynamics}, we describe necessary terminology from
complex dynamics and give a few fundamental lemmas.
Section~\ref{SecPeriodic} contains a proof of the periodic part of
the theorem, and along the way it shows how to interpret the
Mandelbrot set as a parameter space of kneading sequences. The
preperiodic part of the theorem is then proved in
Section~\ref{SecPreperiodic}, using properties of kneading sequences.
In the final Section~\ref{SecHyperbolic}, we show that the methods of
Section~\ref{SecPeriodic} can be used to prove an Orbit Separation
Lemma, which has interesting consequences about hyperbolic components
of the Mandelbrot set. Most of the results and proofs in this paper
work also for ``Multibrot sets'': these are the connectedness loci of
the polynomials $z^d+c$ for $d\geq 2$.

This article is an elaborated version of Chapter~2 of my Ph.D.\
thesis~\cite{DSThesis} at Cornell University, written under the
supervision of John Hubbard and submitted in the summer of 1994. It is
part of a mathematical ping-pong  with John Milnor: it builds at
important places on the paper \cite{GM}; recently Milnor has written
a most beautiful new paper \cite{MiOrbits} investigating external
rays of the Mandelbrot set from the point of view of ``orbit
portraits'', i.e., landing patterns of periodic dynamic rays. I have
not tried to hide how much both I and this paper have profited from
many discussions with him, as will become apparent at many places.
This paper, as well as Milnor's new one, uses certain global counting
arguments to provide estimates, but in different directions. It is a
current project to combine both approaches to provide a new, more
conceptual proof without global counting. Proofs in a similar spirit of
further fundamental properties of the Mandelbrot set can be found in
\cite{DS_Fibers}.

{\sc Acknowledgements.}
It is a pleasure to thank many people who have contributed to
this paper in a variety of ways. I am most grateful to John Hubbard
for so much help and friendship over the years. To John Milnor, I am
deeply indebted for the ping-pong mentioned above, as well as for 
many discussions and a lot of encouragement along the way. Many more
people have shared their ideas and understanding with me and will
recognize their contributions; they include Bodil Branner, Adrien
Douady, Karsten Keller, Eike Lau, Jiaqi Luo, Misha Lyubich, Shizuo
Nakane, Chris Penrose, Carsten Petersen, Johannes Riedl, Mitsu
Shishikura and others.  I am also grateful to the Institute for
Mathematical Sciences in Stony Brook for its inspiring environment and
support. Finally, special thanks go to Katrin Wehrheim for a most
helpful suggestion.

\newsection{Complex Dynamics}
\label{SecDynamics}

In this section, we briefly recall some results and notation from
complex dynamics which will be needed in the sequel. For details, the
notes \cite{MiIntro} by Milnor are recommended and, of course, the
work \cite{Orsay} by Douady and Hubbard which is the source of most
of the results mentioned below.

By affine conjugation, quadratic polynomials can be written uniquely
in the normal form $p_c:z\mapsto z^2+c$ for some complex parameter $c$.
For any such polynomial, the filled-in Julia set is defined as the
set of points $z$ with bounded orbits under iteration. The Julia set
is the boundary of the filled-in Julia set. It is also the set of
points which do not have a neighborhood in which the sequence of
iterates is normal (in the sense of Montel). Julia set and filled-in
Julia set are connected if and only if the only critical point $0$
has bounded orbit; otherwise, these sets coincide and are a Cantor
set. The Mandelbrot set $\M$ is the {\em quadratic connected locus:}
the set of parameters $c$ for which the Julia set is connected.
Julia sets and filled-in Julia sets, as well as the Mandelbrot set,
are compact subsets of the complex plane. The Mandelbrot set is known
to be connected and full (i.e. its complement is connected). 

Douady and Hubbard have shown that Julia sets and the Mandelbrot set
can profitably be studied using external rays: for a compact
connected and full set $K\subset\cz$, the Riemann mapping theorem
supplies a unique conformal isomorphism $\Phi_K$ from the exterior of
$K$ to the exterior of a unique disk $\ovl{D}_R=\{z\in\cz:|z|\leq R\}$
subject to the normalization condition
$\lim_{z\rightarrow\infty}\Phi(z)/z= 1$. The inverse of the Riemann map
allows to transport polar coordinates to the exterior of $K$; images of
radial lines and centered circles are called {\em external rays} and
{\em equipotentials}, respectively. For a point $z\in\cz-K$ with
$\Phi(z)=r e^{2\pi i\theta}$, the number $\theta$ is called the {\em
external angle} and $\log r$ is called the {\em potential} of $z$.
External angles live in $\Circle$; we will always measure them in full
turns, i.e., interpreting $\Circle=\rz/\zz$. Sometimes, it will be
convenient to count the two angles $0$ and $1$ differently and have
external angles live in $[0,1]$. Potentials are parametrized by the
open interval $(\log R,\infty)$. An external ray at angle $\theta$ is
said to {\em land} at a point $z$ if $\lim_{r\searrow \log R}
\Phi_K^{-1}(re^{2\pi i\theta})$ exists and equals $z$. For general
compact connected full sets $K$, not all external rays need to land. By
Carath\'eodory's theorem, local connectivity of $K$ is equivalent to
landing of all the rays, with the landing points depending continuously
on the external angle. 

For all the sets we consider here, it turns out that the conformal
radius $R$ is necessarily equal to $1$. In order to avoid confusion,
we will replace the term ``external ray'' by {\em dynamic ray} or
{\em parameter ray} according to whether it is an external ray of
a filled-in Julia set or of the Mandelbrot set. 

For $c\in\M$, the filled-in Julia set $K_c$ is connected. For
brevity, we will denote the preferred Riemann map by $\phi_c$, rather
than $\Phi_{K_c}$. A classical theorem of B\"ottcher asserts that
this map conjugates the dynamics outside of $K_c$ to the squaring map
outside the closed unit disk: $\phi_c\circ p_c = (\phi_c)^2$.
A dynamic ray is periodic or preperiodic whenever its external angle
is periodic or preperiodic under the doubling map on $\Circle$. The
periodic and preperiodic angles are exactly the rational numbers.
More precisely, a rational angle is periodic iff, when written in
lowest terms, the denominator is odd; if the denominator is even,
then the angle is preperiodic. It is well known
\cite[Section~18]{MiIntro} that dynamic rays of connected filled-in
Julia sets always land whenever their external angles are rational.
The landing points of periodic (resp.\ preperiodic) rays are periodic
(resp.\ preperiodic) points on repelling or parabolic orbits.
Conversely, every repelling or parabolic periodic or preperiodic point
of a connected Julia set is the landing point of one or more rational
dynamic rays; preperiods and periods of all the rays landing at the
same point are equal.

If a quadratic Julia set is a Cantor set, then there still is a
B\"ottcher map $\phi_c$ near infinity conjugating the dynamics to the
squaring map. One can try to extend the domain of definition of the
B\"ottcher map by pulling it back using the conjugation relation.
However, there are problems about choosing the right branch of a
square root needed in the conjugation relation. The absolute value of
the B\"ottcher map is independent of the choices and allows to define
potentials outside of the filled-in Julia set. The set of points at
potentials exceeding the potential of the critical point is simply
connected and the map $\phi_c$ can be defined there uniquely. This
domain includes the critical value. In particular, the external angle
of the critical value is defined uniquely. Douady and Hubbard have
shown that the preferred Riemann map $\Phi_\M$ of the exterior of the
Mandelbrot set is given by $\Phi_\M(c)=\phi_c(c)$. 

For disconnected Julia sets, the map $\phi_c$ defines dynamic rays at
sufficiently large potentials. If a dynamic ray at angle $\theta$ is
defined for potentials greater than $t>0$, then one can pull back by
the dynamics and obtain the dynamic rays at angles $\theta/2$ and
$(\theta+1)/2$ down to potential $t/2$, except if the ray at
angle $\theta$ contains the critical value. In the latter case, the
two pull-back rays will bounce into the critical point and the
pull-back is no longer possible uniquely. This phenomenon has been
studied by Goldberg and Milnor in the appendix of \cite{GM}.
Conversely, a dynamic ray at angle $\theta$ can be extended down to
the potential $t>0$ provided its image ray at angle $2\theta$ can be
extended down to the potential $2t$ and does not contain the critical
value, or if the ray at angle $4\theta$ can be extended down to the
potential $4t$ without containing the critical value or its image,
etc.. The ray can be defined for all potentials in $(0,\infty)$ if
the external angle of the critical value is different from
$2^k\theta$ for all $k=1,2,3,\ldots$. This is the general situation,
and in this case, the dynamic ray is known to land at a unique point
of the Julia set, whether or not the angle $\theta$ is rational. 

We rephrase these facts in a form which we will have many
opportunities to use: if a parameter $c\notin\M$ has external angle
$\theta$, then the dynamic ray at angle $\theta$ for the parameter
$c$ will contain the critical value. If the angle $\theta$ is
periodic, then this ray cannot possibly land: the ray must bounce
into an inverse image of the critical point at a finite positive
potential. The main focus of Sections~\ref{SecPeriodic} and
\ref{SecPreperiodic} will be to transfer the landing properties of
dynamic rays at rational angles into landing properties of parameter
rays at rational angles: as so often in complex dynamics, the general
strategy is ``to plow in the dynamical plane and then to harvest in
parameter space'', as Douady phrased it.

When a periodic ray lands at a periodic point, the periods need not
be equal: it is possible that the period of the ray is a proper
multiple of the period of the point it is landing at. We will
therefore distinguish {\em ray periods} and {\em orbit periods}. If
only one ray lands at every periodic point on the orbit, then both
periods are equal; in general, there is a relation between these
periods and the number of rays landing at each point on the orbit;
see Lemma~\ref{LemRayPermutation}. For our purposes, periodic
orbits will be most interesting if at least two rays land at each of
its points. Such periodic orbits have a distinguished point and two
distinguished dynamic rays landing at this point; these play a
prominent role in all the symbolic descriptions of the Mandelbrot set.
Following the terminology of Milnor~\cite{MiOrbits}, we will call the
distinguished point and rays the {\em characteristic periodic point of
the orbit} and the {\em characteristic rays} (see below), and the
corresponding external angles will be the two {\em characteristic
angles} of the orbit. In Thurston's fundamental preprint \cite{Th}, the
two characteristic rays and their common landing point are the ``minor
leaf'' of a ``lamination''. We will not use or describe his notation
here, but we note that it is very close in spirit to this article.

For our purposes, it will be sufficient to define characteristic
points and rays only for parabolic periodic orbits. 

\begin{definition}[Characteristic Components, Points and Rays]
\label{DefCharacteristic} \lineclear
For a quadratic polynomial with a parabolic orbit, the unique Fatou
component containing the critical value will be called the {\em
characteristic Fatou component}; the only parabolic periodic point on
its boundary will be the {\em characteristic periodic point} of the
parabolic orbit. It is the landing point of at least two dynamic rays,
and the two of them closest to the critical value on either side will
be the {\em characteristic rays}. 
\end{definition}
The fact that every parabolic periodic point is the landing point of at
least two dynamic rays will be shown after
Lemma~\ref{LemHubbardBranch}. Lemma~\ref{LemRayPermutation} will
describe the characteristic rays dynamically.

With hesitation, we use the term ``Misiurewicz point'' for a
parameter $c$ for which the critical point or, equivalently, the
critical value, is (strictly) preperiodic. This terminology has been
introduced long ago, but it is only a very special case of what
Misiurewicz was investigating. In real dynamics, the term is
used in a wider meaning. We have not been successful in finding 
an adequate substitution term and invite the reader for
suggestions.

In this section, we provide two lemmas which are the engine of our
proof: the first one is of analytical nature; it is a slight
generalization of Lemma~B.1 in Goldberg and Milnor~\cite{GM},
guaranteeing stability in the Julia set at repelling (pre)periodic
points. The second lemma will make counting possible by estimating
the number of parabolic parameters with given ray periods.


\begin{lemma}[Stability of Repelling Orbits] 
\label{LemStable} \lineclear 
Suppose that, for some parameter $c_0\in\cz$ (not necessarily in the
Mandelbrot set), there is a repelling periodic point $z_0$ at which
some periodic dynamic ray at angle $\theta$ lands. Then, for
parameters $c$ sufficiently close to $c_0$, the periodic point $z_0$
can be continued analytically as a function $z(c)$ and the dynamic
ray at angle $\theta$ in the dynamic plane of $c$ lands at $z(c)$.
Moreover, the dynamic ray and its landing point form a closed set which
is canonically homeomorphic to $[0,\infty]$ via potentials, and this
parametrized ray depends continuously on the parameter. 

If $z_0$ is repelling and preperiodic, the analogous statement holds
provided that neither the point $z_0$ nor any point on its forward
orbit is the critical point. 
\end{lemma} 

\proof 
We first assume $z_0$ to be a periodic point.
By the implicit function theorem, $z_0$ can be continued analytically
as a function $z(c)$ in a neighborhood of $c_0$; the multiplier
$\lambda(c)$ will also depend analytically on $c$ so that the cycle
is repelling sufficiently close to $c_0$. Let $V$ be such a
neighborhood of $c_0$ and denote the period of $z_0$ by $n$. Then for
every $c\in V$ there exists a local branch $g_c$ of the inverse map
of $p_c^{\circ n}$ fixing $z(c)$. There is a neighborhood $U$ of
$z_0$ such that $g_{c_0}$ maps $\ovl U$ into $U$, and possibly by
shrinking $V$, we may assume that all $g_c$ have the same property
for $c\in V$. Under iteration of $g_c$, any point in $U$ then
converges to $z(c)$. Let $t>0$ be a potential such that, for the
parameter $c_0$, the set $U$ contains all the points of the dynamic
$\theta$-ray at potentials $t$ and below, including the landing point.

Now we distinguish two cases, according to whether or not $c_0\in\M$.
If $c_0\notin\M$, then the external angle of the parameter $c_0$ is
well-defined and different from the finitely many angles $2^k\theta$
for $k=1,2,3,\ldots$ because the dynamic ray at angle $\theta$ lands.
If $V$ is small enough so that all points in $V$ are outside $\M$ and
have their external angles different from all the $2^k\theta$, then
for every $c\in V$, the dynamic ray at angle $\theta$ lands, and the
point at potential $t$ depends analytically on the parameter. It will
therefore be contained in $U$ for sufficiently small perturbations
and thus converge to $z(c)$ under iteration of $g_c$, so the landing
point of the ray is $z(c)$. 

However, if $c_0\in\M$, then we may assume $V$ small enough so that
all its points have potentials less than $t/2$ (with respect to
the potential function of the Mandelbrot set). In the corresponding
dynamic planes, the critical values then have potentials less than
$t/2$, so every dynamic ray exists and depends analytically on
the parameter for potentials greater than $t/2$. By shrinking
$V$, we may then assume that for all $c\in V$, the segment between
potentials $t/2$ and $t$ in the dynamic ray at angle $\theta$ is
contained in $U$. Iterating the map $g_c$, it follows that the
dynamic ray at angle $\theta$ lands at $z(c)$. In both cases, rays
and landing points depend continuously on the parameter, including the
parametrization by potentials. 

The statement about preperiodic points follows by taking inverse
images and is straightforward, except if $z_0$ or any point on its
forward orbit are the critical point. However, if some preperiodic
dynamic ray lands at the critical value, then a small perturbation
may bring the critical value onto this dynamic ray, and the inverse
rays will bounce into the critical point (after that, both branches
will land, and the landing points are two branches of an analytic
function). 
\qed

\begin{lemma}[Counting Parabolic Orbits] 
\label{LemParaCount} \lineclear
For every positive integer $n$, the number of parabolic parameters
in $\cz$ having a parabolic orbit of exact ray period $n$ is at most
half the number of periodic angles in $[0,1]$ having exact period
$n$ under doubling modulo $1$.
\end{lemma}
\proof
We can calculate the exact number of periodic angles. If an angle
$\theta\in[0,1]$ satisfies $2^n\theta\equiv\theta$ modulo $1$, then
we can write $\theta=a/(2^n-1)$ for some integer $a$, and there are
$2^n$ such angles in $[0,1]$. Only a subset of these angles has
exact period $n$: denoting the number of such angles by $s'_n$, we
have $\sum_{k|n} s'_k = 2^n$, which allows to determine the $s'_n$
recursively or via the M\"obius inversion formula. We have $s'_1=2$,
and all the $s'_n$ are easily seen to be even. In the sequel, we will
work with the integers $s_n:=s'_n/2$. The first few terms of the
sequence $(s_n)$, starting with $s_1$, are $1, 1, 3, 6, 15, 27, 63,
\ldots$. The specified number of periodic angles in $[0,1]$ is then
exactly $2s_n$. 

We consider the curve $\{(c,z)\in\cz^2: p_c^{\circ n}(z)=z\}$
consisting of points $z$ which are periodic under $p_c$ with period
dividing $n$. It factors as a product $\prod_{k|n} Q_k(c,z)$
according to exact periods. (The curves $Q_k$ have been shown to be
irreducible by Bousch~\cite{Bousch} and by Lau and
Schleicher~\cite{IntAddr}, a fact we will not use.) For $|c|>2$, the
filled-in Julia set of $p_c$ is a Cantor set containing all the
periodic points. For $|c|>4$, it is easy to verify that points $z$
with $|z|>|c|^{1/2}+1$ escape to $\infty$, and so do points with
$|z|<|c|^{1/2}-1$. Periodic points therefore satisfy
$|z|=|c|^{1/2}(1+o(1))$ as $c\rightarrow\infty$. The multiplier of a
periodic orbit of exact period $n$ is the product of the periodic
points on the orbit multiplied by $2^n$, so it grows like
$|4c|^{n/2}(1+o(1))$. 

For any parameter $c$, the number of points which are fixed under the
$n$-th iterate is obviously equal to $2^n$, counting multiplicities.
These points have exact periods dividing $n$, so the number of
periodic points of exact period $n$ equals $2s_n$ by the same
recursion formula as above. These periodic points are grouped in
orbits, so the number of orbits is $2s_n/n$ (which implies that $2s_n$
is divisible by $n$). For bounded parameters $c$, the periodic points
and thus the multipliers are bounded; since there are $2s_n/n$ orbits,
the multipliers of which are analytic and behave like $|c|^{n/2}$ near
infinity, and since every orbit contains $n$ points, it follows that
sufficiently large multipliers are assumed exactly
$(2s_n/n)(n/2)n=ns_n$ times on $Q_n$ (we do not have to count
multiplicities here because multiple orbits always have multiplier
$+1$). Consider the multiplier map on $Q_n$ which assigns to every
point $(c,z)$ the multiplier $(\partial/\partial z)p_c^{\circ n}(z)$. 
It is a proper map and thus has a mapping degree, so (counting
multiplicities) every multiplier in $\cz$ is assumed equally often,
including the value $+1$. The number of points $(c,z)$ having
multiplier $+1$ therefore equals $ns_n$, counting multiplicities.
Projecting onto the $c$-coordinate and ignoring multiplicities, a
factor $n$ is lost because points on the same orbit project onto the
same parameter, and we obtain an upper bound of $s_n$ for the number of
parameters. (In fact, it is not too hard to show at this point that
$s_n$ provides an exact count \cite{MiOrbits}. We will show this in
Corollary~\ref{CorParaCountExact} by a global counting argument.)

Consider a parabolic orbit of exact period $k$ and multiplier
$\mu=e^{2\pi ip/q}$ with $(p,q)=1$. Then the exact ray period is
$qk=:n$, and $qk$ is also the smallest period such that, when
interpreting the orbit as an orbit of this period, the multiplier
becomes $+1$. Therefore, the periodic points on this orbit are on
$Q_{n}$, and the number of parabolic parameters having exact ray
period $n$ therefore is at most $s_n$.
\qed

A more detailed account of such counting arguments can be found in
Section~5 of Milnor~\cite{MiOrbits}.

\bigskip
The following standard lemma is folklore and at the base of every
description of quadratic iteration theory. Our proof follows
Milnor~\cite{MiOrbits}; compare also
Thurston~\cite[Theorem~II.5.3 case b) i) a)]{Th}. We do not assume
the Julia set to have any particular property; it need not even be
connected.

\goodbreak
\begin{lemma}[Permutation of Rays]\nobreak
\label{LemRayPermutation} \lineclear
If more than two periodic rays land at a periodic point, or if the
orbit period is different from the ray period, then the first return
map of the point permutes the rays transitively.
\end{lemma}
\proof
Denote the orbit period by $k$ and the ray period by $n$. Since a
periodic orbit has periodic rays landing only if the orbit is
repelling or parabolic, the first return map of any of its periodic
points is a local homeomorphism and permutes the rays landing there
in such a way that their periods are all equal, and the number of
rays landing at each point of the orbit is a constant $s$, say. If
$s=1$, then orbit period and ray period are equal. If $s=2$, then
either ray and orbit periods are equal, or the first return map of
any point has no choice but to transitively permute the two rays
landing at this point. We may hence assume $s\geq 3$. Then the $s$
rays landing at any one of these periodic points separate the dynamic
plane into $s$ sectors. Every sector is bounded by two dynamic rays,
so it has associated a width: the external angles of the two rays cut
$\Circle$ into two open intervals, exactly one of which does not
contain external angles of rays landing at the same point. The {\em
width} of the sector will be the length of this interval (normalized
so that the total length of $\Circle$ is $1$). 

Since the dynamics of the first return map is a local homeomorphism
near the periodic point, every sector is periodic, and so is the
sequence of the corresponding widths. More precisely, we will justify
the following observations below: if a sector does not contain the
critical point, then it maps homeomorphically onto its image sector
(based at the image of its landing point), and the width of the
sector doubles. However, if the sector does contain the critical
point, then the sector maps in a two-to-one fashion onto the image
sector, and it covers the remaining dynamic plane once. In this case,
the width of the sector will decrease under this mapping, and the
image sector contains the critical value. To justify these statements,
first note that the rays bounding any sector are mapped to the rays
bounding the image sector. Looking at external angles within the
sector, it follows that either the sector maps forward
homeomorphically, or it covers the entire complex plane once and the
image sector twice. The latter must happen for the sector containing
the critical point. Since all the sectors at any periodic point
combined exactly cover the complex plane twice when mapped forward,
all the other sectors must map homeomorphically onto the image
sectors. We also see that among all the sectors based at any point,
the sector containing the critical point must have width greater than
$1/2$, and all the other sectors then have widths less than $1/2$
(the critical point cannot be on a sector boundary: if it is on a
periodic dynamic ray, then this ray cannot land, and if it is on a
periodic point, then this point is superattracting). The width of any
sector doubles under the map if it does not contain the critical point;
since the sum of the widths of all the sectors based at any point
is $1$, the width of the critical sector must decrease.

For each orbit of sectors, there must be at least one sector with
minimal width. It must contain the critical value (or it would be the
image of a sector with half the width), and it cannot contain the
critical point (or its image sector would have smaller width).
Therefore, all the shortest sectors of the various cycles of sectors
must be bounded by pairs of rays separating the critical point from
the critical value, and these sectors are all nested. Among them,
there is one innermost sector $S_1$ based at some point $z_1$ of the
periodic orbit. This sector $S_1$ cannot contain another point from the
orbit of $z_1$: if there was such a point $z'$, there would have to be
a sector based at $z'$ which was shorter than all the shortest sectors
at points on the orbit of $z_1$, and this is obviously absurd. 

\hide{
Now we show that $S_1$ does not contain another point from the orbit
of $z_1$. If it does, let $z'$ be such a point and let $S'$ be a
sector based at $z'$ which does not contain $z_1$; it cannot contain
the critical point, either. By minimality of $S_1$, the sector $S'$
cannot contain the critical value. There is a shortest sector on the
orbit of $S'$. It does contain the critical value and thus the sector
$S_1$ (it might be the sector $S_1$). But then it also contains $S'$
and must then have larger width than $S'$. This contradiction shows
that $S_1$ cannot contain any point in the orbit of $z_1$.
}

If there is an orbit of sectors not involving $S_1$, then any
shortest sector on this orbit must contain the critical value and
thus $S_1$, but it cannot contain the critical point. This sector
must then contain all sectors at $z_1$ except the one containing the
critical point. The representative of this orbit of sectors at $z_1$
must then be the unique sector containing the critical point. Any
cycle of sectors has then only two choices for its representative at
$z_1$: the sector containing the critical point or the critical value.
If there is more than one cycle, then it follows that there are just
two cycles, and each of them has exactly one representative at each
point of the periodic orbit, so $s=2$ in contradiction to our
assumption.
\qed

\remark
This lemma is at the heart of the general definition of
characteristic rays: the main part of the proof works when at least
two rays land at each of the periodic points, and it shows that there
is a unique sector of minimal width containing the critical value.
The rays bounding this sector are called the characteristic rays.
For the special case of a parabolic orbit, this definition agrees
with the one we have given above.

\newsection {Periodic Rays}
\label{SecPeriodic}

In this section, we will be concerned with parameter rays at periodic
angles. The proof of the following weak form of the theorem is due to
Goldberg, Milnor, Douady, and Hubbard; see \cite[Theorem~C.7]{GM}.

\begin{proposition}[Periodic Parameter Rays Land]
\label{PropPeriodicRaysLand} \lineclear
Every parameter ray at a periodic angle $\theta$ lands at a parabolic
parameter $c_0$. In the dynamic plane of\/ $c_0$, the dynamic ray at
angle\/ $\theta$ lands at the parabolic orbit. 
\end{proposition} 
\proof
Let $c_0$ be a parameter in the limit set of the parameter ray at
angle $\theta$ and let $n$ be the exact period of $\theta$. In the
dynamic plane of $c_0$, the dynamic ray at angle $\theta$ must land
at a repelling or parabolic periodic point $z$ of ray period $n$; see
\cite[Theorem~18.1]{MiIntro}. If $z$ was repelling,
Lemma~\ref{LemStable} would imply that for parameters $c$ sufficiently
close to $c_0$, the dynamic ray at angle $\theta$ in the dynamic plane
of $c$ would land at a repelling periodic point $z(c)$, so it could not
bounce off any precritical point. However, when $c$ is on the parameter
ray at angle $\theta$, then the dynamic ray at angle $\theta$ must
bounce off some precritical point, even infinitely often. 

Therefore, $c_0$ is parabolic, and within its dynamics, the dynamic
ray at angle $\theta$ lands at the parabolic orbit. Since limit sets
are connected but parabolic parameters of given ray period form a
finite set by Lemma~\ref{LemParaCount}, the parameter ray at angle
$\theta$ lands and the statements follow. 
\qed

This proves half of the first assertion in Theorem~\ref{ThmRatRays}.
The remainder of the first and the second assertion will be shown in
several steps. We want to show that at a parabolic parameter $c_0$,
those two parameter rays land which have the same external angles as
the two characteristic rays of the critical value Fatou component,
and no other rational ray lands there. The first statement is usually
shown using Ecalle cylinders. It turns out that it is much easier to
show that some ray does {\em no\/t} land at a given point, rather
than to show where it does land. The idea in this paper will be to
exclude all the wrong rays from landing at given parabolic
parameters, using partitions in the dynamic and parameter planes.
Using that the rays must land somewhere, a global counting argument
will then prove the theorem. 

Let $c$ be a parabolic parameter and let $\Theta_c$ be the set of
periodic angles $\theta$ such that the parameter ray at angle
$\theta$ lands at $c$. A priori, it might be empty. We will prove the
following two results later in this section.

\begin{proposition}[Necessary Condition]
\label{PropNecessary} \lineclear
If an angle $\theta$ is in $\Theta_c$, then the dynamic ray at angle
$\theta$ lands at the characteristic point of the parabolic orbit in
the dynamic plane of $c$. 
\end{proposition}

\begin{proposition}[At Most Two Rays]
\label{PropAtMostTwoRays} \lineclear
If\/ $\Theta_c$ contains more than one angle, then it consists of
exactly those two angles which are the characteristic angles of the
parabolic orbit in the dynamic plane of $c$. 
\end{proposition}

These two propositions allow to prove the half of the theorem dealing
with periodic rays; we will deal with the preperiodic half in the
next section.
\proofof{Theorem~\ref{ThmRatRays} (periodic case)}
By Lemma~\ref{LemParaCount}, the number of parabolic parameters of
any given ray period is at most half the number of parameter rays at
periodic angles of the same period. Since every ray lands at such a
parabolic parameter by Proposition~\ref{PropPeriodicRaysLand}, and at
most two rays may land at any such point by
Proposition~\ref{PropAtMostTwoRays}, it follows that exactly two rays
land at every parabolic point, and Proposition~\ref{PropAtMostTwoRays}
says which ones these are. It also follows that the number of
parabolic parameters of any given period is largest possible as
allowed by Lemma~\ref{LemParaCount}.
\qed
\remark
Since we will complete the proof of Proposition~\ref{PropAtMostTwoRays}
by induction on the period, using Theorem~\ref{ThmRatRays} for lower
periods, it is important to note that in order to prove the Theorem
for any period, it suffices to know Proposition~\ref{PropAtMostTwoRays}
for the same period.

\begin{corollary}[Counting Parabolic Orbits Exactly]
\label{CorParaCountExact} \lineclear
Let $s_k$ be the number of parameters having a parabolic orbit of
exact ray period $k$. These numbers satisfy the recursive relation
$\sum_{k|n} s_k = 2^{n-1}$, which determines them uniquely.
\qedd
\end{corollary}

\begin{figure}[htbp]
\centerline{\psfig{figure=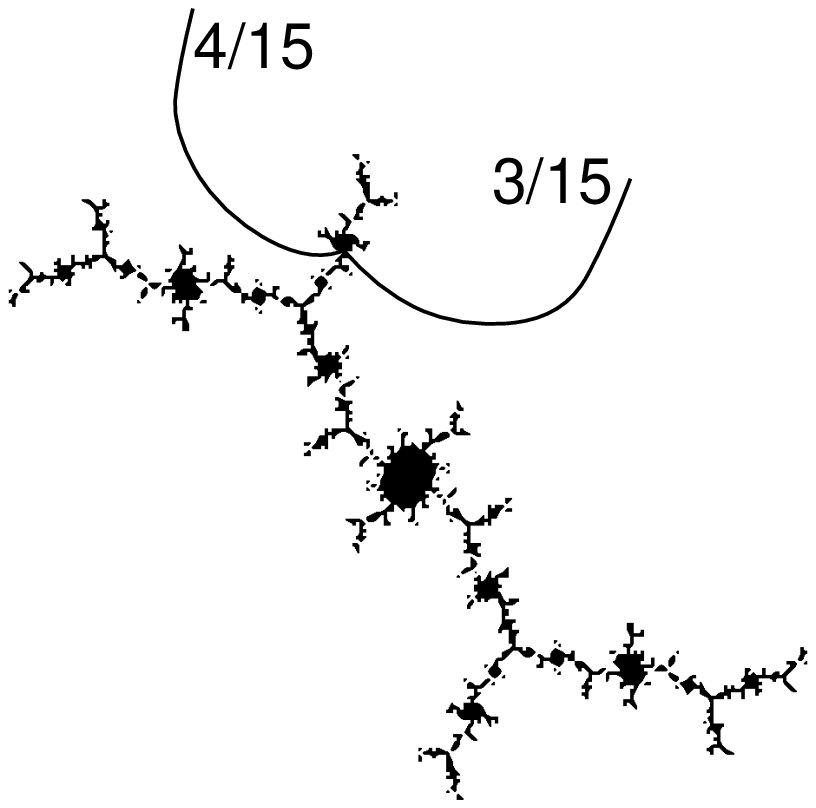,height=55mm} \hfil
            \psfig{figure=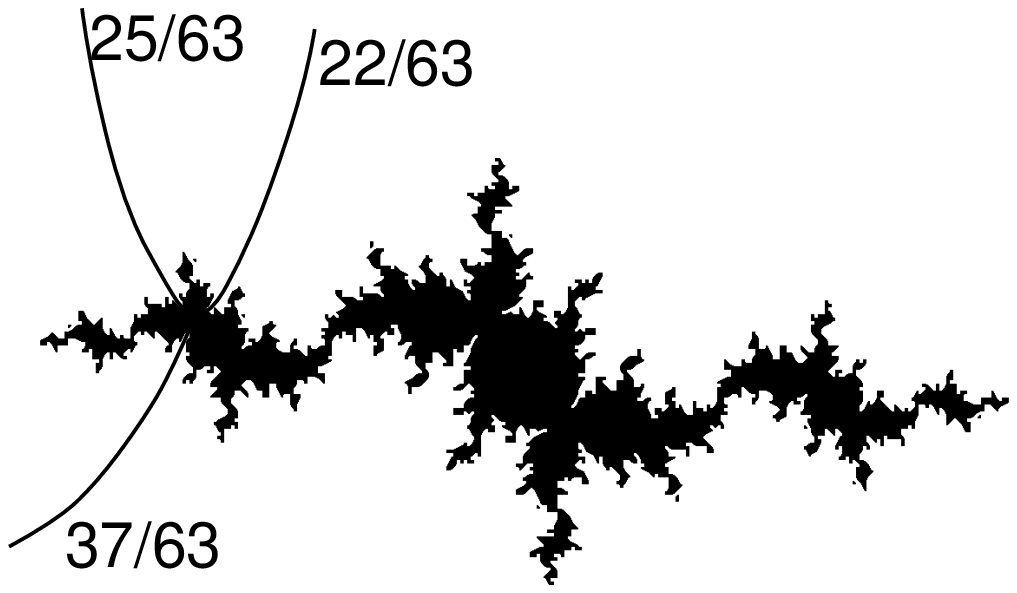,height=55mm}}
\centerline{\parbox{\captionwidth}{ 
\caption{\sl Illustration of the theorem in the periodic case. The
polynomials at the landing points of the parameter rays at angles
$3/15$ and $4/15$ (left) and at angles $22/63$ and $25/63$ (right)
are shown. In both pictures, the rays landing at the characteristic
points are drawn. For the corresponding parameter rays, see
Figure~\protect\ref{FigMandelRays}. }
\label{FigPeriodic}
} }
\end{figure}

It remains to prove the two propositions. In both of them, we have to
exclude that certain rays land at given parabolic parameters. We do
that using appropriate partitions: first in the dynamic plane, then
in parameter space. We start by discussing the topology of parabolic
quadratic Julia sets and define a variant of the Hubbard tree on
them. Hubbard trees have been introduced by Douady and Hubbard in
\cite{Orsay} for postcritically finite polynomials. We will be
interested in combinatorial statements about combinatorially described
Julia sets, so these results could be derived in purely combinatorial
terms. However, it will be more convenient to use topological
properties of the Julia sets in the parabolic case, in particular that
they are pathwise connected (which follows from local connectivity).
This was originally proved by Douady and Hubbard~\cite{Orsay}; proofs
can also be found in Carleson and Gamelin~\cite{CG} and in Tan and
Yin~\cite{TY}. 

In a quadratic polynomial with a parabolic orbit, let $z$ be any point
within the filled-in Julia set and let $U$ be a bounded Fatou
component. We then define a {\em projection} of $z$ onto $\ovl U$ as
follows: If $z\in\ovl U$, then the projection of $z$ onto $\ovl U$
is $z$ itself. Otherwise, consider any path within the filled-in
Julia set connecting $z$ to an interior point of $U$ (such a path
exists because the filled-in Julia set is pathwise connected); then
the projection of $z$ onto $\ovl U$ is the first point where this path
intersects $\partial U$. There may be many such paths, but the
projection is still well-defined: take any two paths from $z$ to the
interior of $U$ and connect their endpoints within $U$. If the paths
are different, they will bound some subset of $\cz$, which must be in
$K$ because $K$ is full. If the paths reach $\partial U$ in different
points, then these paths enclose part of the boundary of $U$, but the
boundary of any Fatou component is always in the boundary of the
filled-in Julia set. This contradiction shows that the projection is
well-defined. Every parabolic periodic point is on the boundary of at
least one periodic Fatou component, so the projection in this case is
just the identity.

\begin{lemma}[Projection Onto Periodic Fatou Components]
\label{LemProject} \lineclear 
In a quadratic polynomial with a parabolic orbit, the projections of
all the parabolic periodic points onto the Fatou component containing
the critical value take images in the same point, which is the
characteristic point of the parabolic periodic orbit. Projections of
the parabolic periodic points onto any other bounded Fatou component
take images in at most two boundary points, which are periodic or
preperiodic points on the parabolic orbit.
\end{lemma}

\proof
Let $n$ be the period of the periodic Fatou components and number them
$U_0, U_1, \ldots,\brkOK U_{n-1},\brkOK U_n=U_0$ in the order of the
dynamics, so 
that $U_0=U_n$ contains the critical point. Let $a_k$ be the number of
different images that the projections of all the parabolic periodic
points onto $\ovl U_k$ have, for $k=0,1,\ldots, n$ (with $a_0=a_n$).
We first show that $a_{k+1}\geq a_k$ for $k=1,2,\ldots, n-1$. 

Let $z$ be a parabolic periodic point and let $\pi(z)$ be its
projection onto $\ovl U_k$. We claim that $p(\pi(z))$ is the
projection onto $\ovl U_{k+1}$ of either $p(z)$ or the parabolic point
on the boundary of the Fatou component containing the critical value
(i.e.\ the characteristic point on the parabolic orbit). Indeed, if the
path between $z$ and $\pi(z)$ maps forward homeomorphically under $p$,
then $\pi(p(z))=p(\pi(z))$. If it does not, then the path must
intersect the component containing the critical point, and
$\pi(z)=\pi(0)$. But then $p(\pi(z))$ is the projection of the
characteristic point on the parabolic orbit. Therefore, for
$k\in\{0,1,2,\ldots, n-1\}$, all the $a_k$ image points of the
projections of parabolic periodic points onto $\ovl U_k$ will be mapped
under $p$ to image points of the projection onto $\ovl U_{k+1}$. Since
for $k\neq 0$, the polynomial $p$ maps $\ovl U_k$ homeomorphically onto
$\ovl U_{k+1}$, we get $a_n\geq a_{n-1}\geq \ldots a_2\geq a_1$.
Similarly, since $p$ maps $\ovl U_0$ in a two-to-one fashion onto $\ovl
U_1$, we have $a_1\geq a_0/2$. 

Now we connect the parabolic periodic points by a tree: first, there
is a path between the critical point and the critical value, and all
the other parabolic periodic points which are not on this path can be
connected, one by one, to the subtree which has been constructed thus
far. We can require that every path which we are adding intersects the
boundary of any bounded Fatou component in the least number of points
(at most two). After finitely many steps, all the parabolic periodic
points are connected by a finite tree, and all the endpoints of this
tree are parabolic periodic points. It is not hard to check that this
tree intersects the boundary of any bounded Fatou component exactly in
the image points of the projections of the parabolic periodic points.

We now claim that there is a periodic Fatou component whose closure
does not disconnect the tree. Indeed, any component which does
disconnect the tree has at least one parabolic periodic point and
thus at least one periodic Fatou component in each connected
component of the complement. Pick one connected component, and within
it pick a periodic Fatou component that is  ``closest'' to the removed
one (in the sense that the path between these two components does not
contain further periodic Fatou components). Remove this component and
continue; this process can be continued until we arrive at a component
which does not disconnect the tree. Let $U_{\tilde k}$ be such a
component.

It follows that $a_{\tilde k}\leq 2$: all the parabolic periodic points
which are not on the boundary of $U_{\tilde k}$ must project to the
same boundary point, say $b$, which may or may not be the parabolic
periodic point on the boundary of $U_{\tilde k}$. We will now show 
that it will be, so $a_{\tilde k}=1$. 

Since $b$ is the image of a projection onto $\ovl U_{\tilde k}$ of
some parabolic periodic point which is not in $\ovl U_{\tilde k}$, it
follows from the argument above that $p(b)$ is the image of a
projection onto $\ovl U_{\tilde k+1}$ of some parabolic periodic point
which is not in $\ovl U_{\tilde k+1}$. But since $b$ is the only
boundary point of $\ovl U_{\tilde k}$ with this property, it follows
that $b$ is fixed under the first return map of this Fatou component.
The point $b$ must then be the unique parabolic periodic point on the
boundary of $\ovl U_{\tilde k}$, and we have $a_{\tilde k}=1$. 

Since $a_1\leq a_2\leq a_n\leq 2a_1$, it follows in particular that
$a_1=1$ and all $a_k\leq 2$, and all projections onto the Fatou
component containing the critical value take values in the same point,
which is the characteristic point of the parabolic orbit. The remaining 
claims follow.
\qed

\remark
The tree just constructed is similar to the {\em Hubbard tree}
introduced in \cite{Orsay} for postcritically finite polynomials. An
important difference is that our tree does not connect the critical
orbit. Moreover, Hubbard trees in \cite{Orsay} are specified uniquely,
while our trees still involve the choice of how to traverse bounded
Fatou components. We will suggest a preferred tree below. 

However, some properties are independent of the choice of the tree.
Assume that two simple curves $\gamma_1$ and $\gamma_2$ within the
filled-in Julia set connect the same two points $z_1$ and $z_2$, such
that a point $w$ is on one of the curves but not on the other. Then $w$
is on the closure of a bounded Fatou component because the region which
is enclosed by the two curves must be in the filled-in Julia set. The
tree intersects the boundary of any bounded Fatou component in at most
two points which are projection images and thus well-defined.
Therefore, the choice for the curves and thus for the tree is only in
the interior of bounded Fatou components.

For any point $w$ in the Julia set (not in a bounded Fatou component),
it follows that the number of branches of the tree (i.e.\ the number of
components the point disconnects the tree into) is independent of the
choice of the tree. Similarly, the number of branches is ``almost''
non-decreasing under the dynamics, so that $p(w)$ has at least as many
branches as $w$: all the different branches at $w$ will yield different
branches at $p(w)$, except if $w$ is on the boundary of the Fatou
component containing the critical point. At such boundary points, only
the branch leading into the critical Fatou component can get lost.
(However, it does happen that $p(w)$ has extra branches in the tree.)
It follows that the characteristic point on the parabolic orbit has at
most one branch on any tree, and all the other parabolic periodic
points can have up to two branches. 

A branch point of a tree is a point $w$ which disconnects the tree into
at least three complementary components. 

\begin{lemma}[Branch Points of Tree]
\label{LemHubbardBranch} \lineclear 
Branch points of the tree between parabolic periodic points are
periodic or preperiodic points on repelling orbits.
\end{lemma}
\proof
Branch points are never on the parabolic orbit, as we have just seen.
Therefore, the image of a branch point is always a branch point with
at least as many branches. Since there are only finitely many branch
points, every branch point is periodic or preperiodic and hence on a
repelling orbit.
\qed

We can now proceed to select a preferred tree, which we will call a
{\em (parabolic) Hubbard tree}. We only have to specify how it will
traverse bounded Fatou components. In fact, since every bounded Fatou
component will eventually map homeomorphically onto the critical Fatou
component, we only have to specify how the tree has to traverse this
component; for the remaining components, we can pull back.

Let $U$ be the critical Fatou component and let $w$ be the parabolic
periodic point on its boundary. First we want to connect the critical
point in $U$ to $w$ by a simple curve which is forward invariant under
the dynamics. We will use {\em Fatou coordinates} for the attracting
petal of the dynamics \cite[Section~7]{MiIntro}. In these coordinates,
the dynamics is simply addition of $+1$, and our curve will just be a
horizontal straight line connecting the critical orbit. This curve can
be extended up to the critical point. The other point on the boundary
of the critical Fatou component which we have to connect is $-w$, and
we use the symmetric curve. With this choice, we have specified a
preferred tree which is invariant under the dynamics, except that the
image of the tree connects the characteristic periodic point on the
parabolic orbit to the critical value. Removing this curve segment from
the image tree, we obtain the same tree as before.

It is well known that, if a repelling or parabolic periodic point
disconnects the Julia set into several parts, then this point is the
landing point of as many dynamic rays as it disconnects the Julia set
into. This is not hard to see once it is known that at least one ray
lands. It follows that any branch point of the Hubbard tree has dynamic
rays landing between any two branches; any periodic point on the
interior of the tree is the landing point of at least two dynamic rays
separating the tree. It now follows that the characteristic point on
the parabolic orbit, and thus every parabolic periodic point, is the
landing point of at least two dynamic rays. The two characteristic rays
of the parabolic orbit are the two rays landing at the characteristic
point of the orbit and closest possible to the critical value on either
side. A different description of the characteristic rays has been given
in Lemma~\ref{LemRayPermutation}.

\begin{lemma}[Orbit Separation Lemma]
\label{LemOrbitSeparation} \lineclear
Any two different parabolic periodic points of a quadratic polynomial 
can be separated by two (pre)periodic dynamic rays landing at a
common repelling (pre)periodic point.
\end{lemma} 
\proof
It suffices to prove the lemma when one of the two parabolic periodic
points is the characteristic point of the orbit; this is also the
only case we will need in this section. Let $z$ be this
characteristic point and let $z'$ be a different parabolic periodic
point. Consider the tree of the polynomial as constructed above. It
contains a unique path connecting $z$ and $z'$. We may assume that
this path does not traverse a periodic Fatou component except at its
ends; if it does, we replace $z'$ by the parabolic periodic point on
that Fatou component. Similarly, we may assume that the path does not
traverse another parabolic periodic point. If the path from $z$ to
$z'$ contains a branch point of the tree, then by
Lemma~\ref{LemHubbardBranch},  this branch point is periodic or
preperiodic and repelling, and it is therefore the landing point of
rational dynamic rays separating the parabolic orbit as claimed.

If the Hubbard tree does not have a branch point between $z$ and
$z'$, then it takes a finite number $k$ of iterations to map $z'$ for
the first time onto $z$. Denoting the path from $z$ to $z'$ by
$\gamma$, then the $k$-th iterate of $\gamma$ must traverse itself
and possibly more in an orientation reversing way: denoting the $k$-th
image of $z$ by $z''$, then the image curve connects $z$ and $z''$; it
must start along the end of $\gamma$ because $z$ is an endpoint of the
Hubbard tree, and it cannot branch off because we had assumed no branch
point of the tree to be on $\gamma$. There must be a unique point
$z_f$ in the interior of $\gamma$ which is fixed under the $k$-th
image of $\gamma$. This point is a repelling periodic point, and it is
the landing point of two dynamic rays with the desired separation
properties.
\qed

Now we can prove Proposition~\ref{PropNecessary}.
\proofof{Proposition~\ref{PropNecessary}}
In Proposition~\ref{PropPeriodicRaysLand}, we have shown that
$\Theta_c$ can contain only angles $\theta$ of dynamic rays landing
at the parabolic cycle in the dynamic plane of $c$. By the Orbit
Separation Lemma~\ref{LemOrbitSeparation}, all the rays not landing
at the characteristic point of the parabolic orbit are separated from
the critical value by a partition formed by two dynamic rays landing
at a common repelling (pre)periodic point. This partition is stable
in a neighborhood in parameter space by Lemma~\ref{LemStable}. But
the parameter $c$ being a limit point of the parameter ray at angle
$\theta$ means that, for parameters arbitrarily close to $c$, the
critical value is on the dynamic ray at angle $\theta$. 
\qed

The set $\Theta_c$ of external angles of the parabolic parameter $c$
can thus contain only such periodic angles which are external angles
of the characteristic periodic point of the parabolic orbit in the
dynamical plane of $c$. If there are more than two such angles, we
want to exclude all those which are not characteristic. This is
evidently impossible by a partition argument in the dynamic plane. 
In order to prove Proposition~\ref{PropAtMostTwoRays}, we will use a
partition of parameter space; for that, we have to look more closely
at parameter space and incorporate some symbolic dynamics using
kneading sequences. The partition of parameter space according to
kneading sequences and, more geometrically, into internal addresses
is of interest in its own right (see below) and has been investigated
by Lau and Schleicher~\cite{IntAddr}; related ideas can be found in
Thurston~\cite{Th}, Penrose~\cite{Pe} and \cite{Pe2} and
in a series of papers by Bandt and Keller (see~\cite{Ke1},
\cite{KeHabil} and the references therein). This partition will also
be helpful in the next section, establishing landing properties of
preperiodic parameter rays.

\begin{definition}[Kneading Sequence]
\label{DefKneadSeq} \lineclear
To an angle $\theta\in\Circle$, we associate its {\em kneading
sequence} as follows: divide\/ $\Circle$ into two parts at\/
$\theta/2$ and\/ $(\theta+1)/2$ (the two inverse images of\/ $\theta$
under angle doubling); the open part containing the angle $0$ is
labeled $\0$, the other open part is labeled $\1$ and the boundary
gets the label $\*$. The kneading sequence of the angle $\theta$ is
the sequence of labels corresponding to the angles $\theta$,
$2\theta$, $4\theta$, $8\theta$, \ldots.
\end{definition}

\begin{figure}[htbp]
\centerline{\psfig{figure=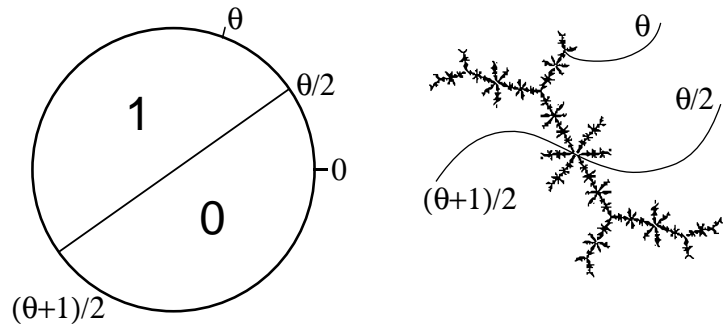,height=45mm}}
\centerline{\parbox{\captionwidth}{ \vspace{0mm} 
\caption{\sl Left: the partition used in the definition of the
kneading sequence. Right: a corresponding partition of the dynamic
plane by dynamic rays, shown here for the example of a Misiurewicz
polynomial. }
\label{FigKneading}
} }
\end{figure}

It is easy to check that, for $\theta\neq 0$, the first position
always equals $\1$. If $\theta$ is periodic of period $n$, then its
kneading sequence obviously has the same property and the symbol $\*$
appears exactly once within this period (at the last position). The
symbol $\*$ occurs only for periodic angles. However, it may happen
that an irrational angle has a periodic kneading sequence
(see e.g.\ \cite{IntAddr}). As the angle $\theta$ varies, the entry
of the kneading sequence at any position $n$ changes exactly at those
values of $\theta$ for which $2^{n-1}\theta$ is on the boundary of
the partition, i.e.\ where the kneading sequence has the entry $\*$.
This happens if and only if the angle $\theta$ is periodic, and its
exact period is $n$ or divides $n$. 

Another useful property which will be needed in
Section~\ref{SecPreperiodic} is that the pointwise limits 
$\K_-(\theta):=\lim_{\theta'\nearrow\theta}\K(\theta')$ and 
$\K_+(\theta):=\lim_{\theta'\searrow\theta}\K(\theta')$ exist for every
$\theta$. If $\theta$ is periodic, then $\K_\pm(\theta)$ is also
periodic with the same period (but its exact period may be smaller).
Both limiting kneading sequences coincide with $\K(\theta)$ everywhere,
except that all the $\*$-symbols are replaced by $\0$ in one of the
two sequences and by $\1$ in the other. The reason is simple: if
$\theta'$ is very close to $\theta$, then the orbits under doubling, as
well as the partitions in the kneading sequences are close to each
other, and any symbol $\0$ or $\1$ at any finite position will be
unchanged provided $\theta'$ is close enough to $\theta$. However, if
the period of $\theta$ is $n$ so that $2^{n-1}\theta$ is on the
boundary of the partition in the kneading sequence, then
$2^{n-1}\theta'$ will barely miss the boundary in its own partition,
and the $\*$ will turn into a $\0$ or $\1$. As long as the orbit of
$\theta'$ is close to the orbit of $\theta$, all the symbols $\*$ will
be replaced by the same symbol.

\begin{figure}[tbh]
\centerline{\psfig{figure=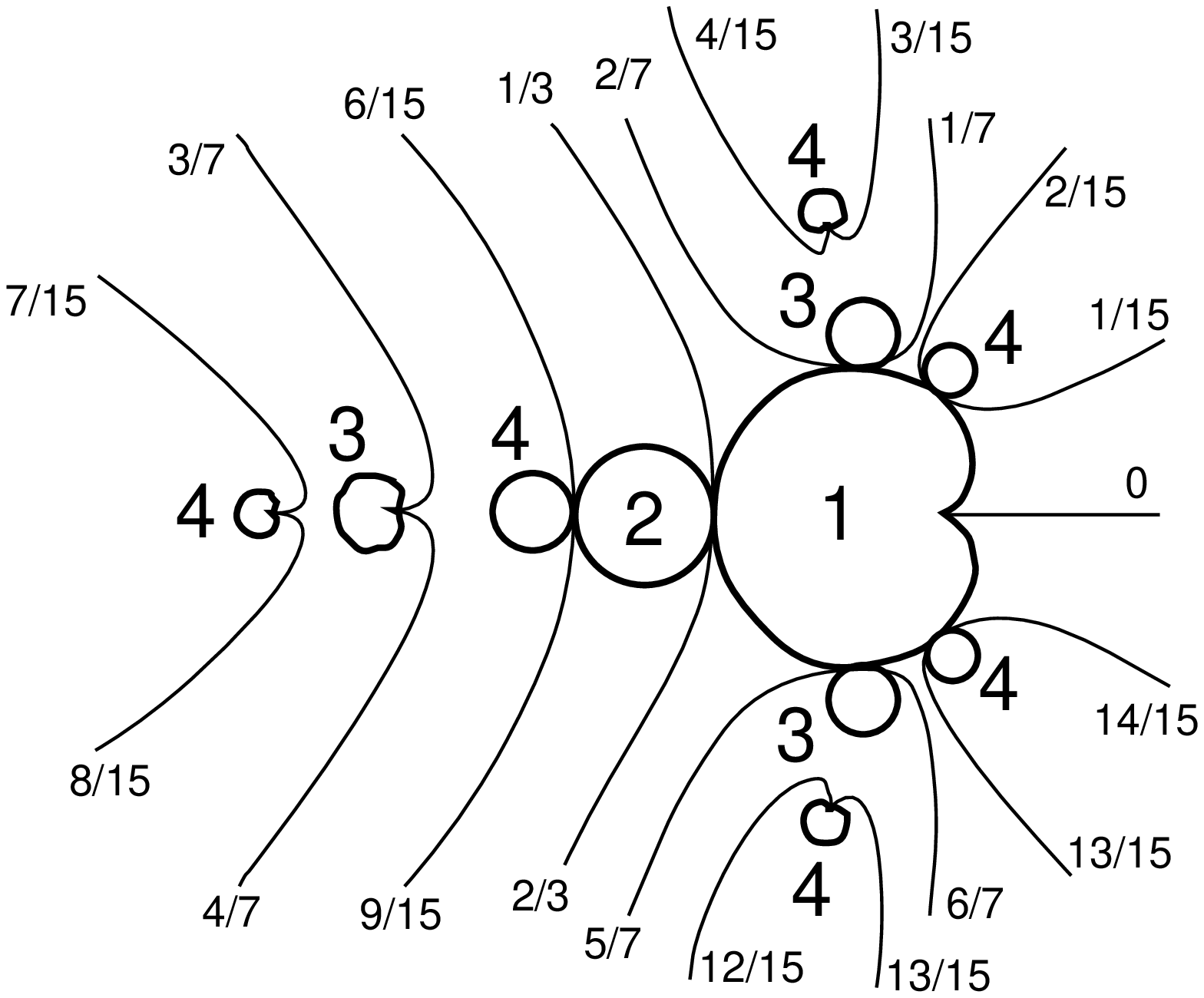,height=70mm} \hfil
            \psfig{figure=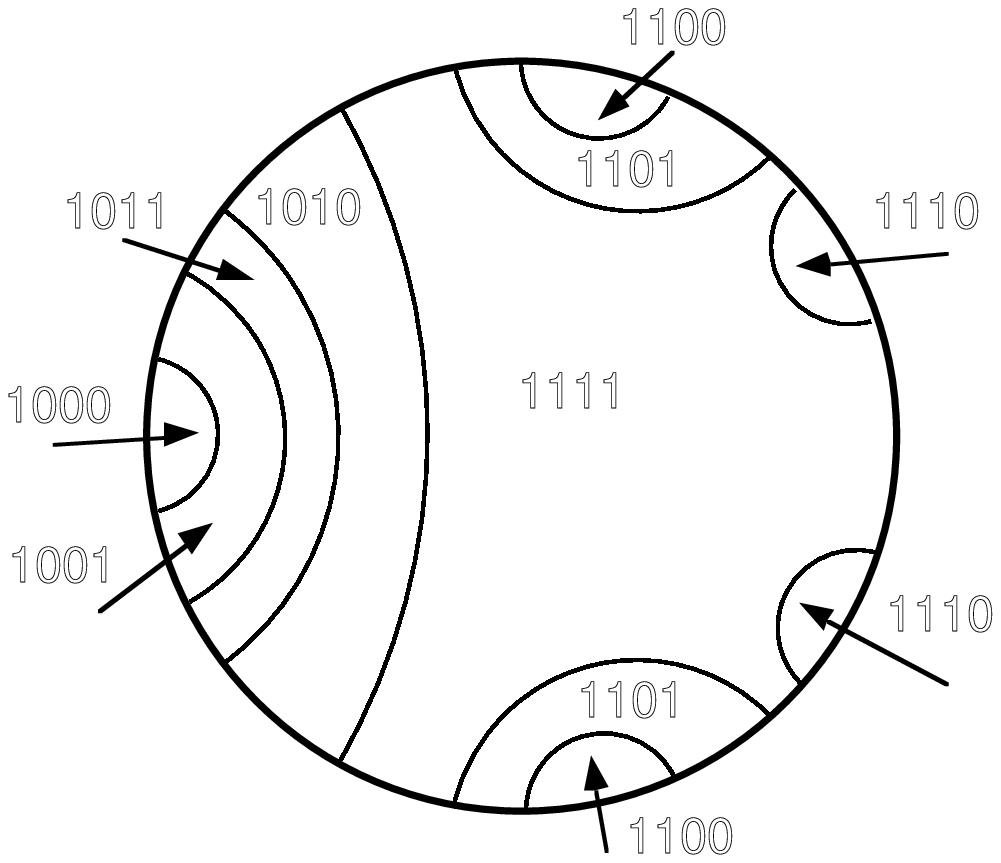,height=70mm}}
\centerline{\parbox{\captionwidth}{ \vspace{3mm} 
\caption[The partition $\Part{n}$ used in the proof of
Theorem~\protect\ref{ThmRatRays}, for $n=4$] {\sl Left: the partition
$\Part{n}$ used in the proof of
Proposition~\protect\ref{PropAtMostTwoRays}, for $n=4$. The
corresponding hyperbolic components are drawn in for clarity and do
not form part of the partition. Right: a corresponding symbolic
picture, showing how the partition yields a parameter space of initial
segments of kneading sequences. The same pairs of rays are drawn in as
on the left hand side, but the angles are unlabeled for lack of
space. }
\label{FigParaPart} } }
\end{figure}

Proposition~\ref{PropPeriodicRaysLand} asserts in particular that all
the periodic parameter rays landing at the same parameter have equal
period. All the rays of period at most $n-1$ divide the plane into
finitely many pieces. We denote this partition by $\Part{n-1}$; it is
illustrated in Figure~\ref{FigParaPart}. Parabolic parameters of ray
period $n$ and parameter rays of period $n$ have no point in common
with the boundary of this partition.

\begin{lemma}[Kneading Sequences in the Partition]
\label{LemKneadingPartition} \lineclear
Fix any period $n\geq 1$ and suppose that all the parameter rays of
periods at most $n-1$ land in pairs. Then all parameter rays in any
connected component of\/ $\Part{n-1}$ have the property that the first
$n-1$ entries in their kneading sequences coincide and do not contain
the symbol $\*$. In particular, rays of period $n$ with different
kneading sequences do not land at the same parameter.
\end{lemma}
\proof
The first statement is trivial for two periodic rays which do not have
rays of lower periods between them, i.e., for rays from the same
``access to infinity'' of the connected component in $\Part{n-1}$:
the first $n-1$ entries in the kneading sequences are stable for
angles within every such access. The claim is interesting only for a
connected component with several ``accesses to infinity''.

The hypothesis of the theorem asserts that parameter rays of periods
up to $n-1$ land in pairs. Therefore, whenever two rays at angles
$\theta_1,\theta_2$ are in the same connected component of
$\Part{n-1}$, the parameter rays of any period $k\leq n-1$ on either
side (in $\Circle$) between these two angles must land in pairs. The
number of such rays is thus even, and the $k$-th entry in the kneading
sequence changes an even number of times between $\0$ and
$\1$. 
\qed

\remark
This lemma allows to interpret $\Part{n}$ as a parameter space of
initial segments of kneading sequences. In Figure~\ref{FigParaPart},
the partition is indicated for $n=4$, together with the initial four
symbols of the kneading sequence. The entire parameter space may thus
be described by kneading sequences, as noted above. To any parameter
$c\in\cz$, we may associate a kneading sequence as follows: it is a
one-sided infinite sequence of symbols, and the $k$-th entry is $\1$
if and only if the parameter is separated from the origin by an even
number of parameter ray pairs of periods $k$ or dividing $k$; if the
number of such ray pairs is odd, then the entry is $\0$, and if the
parameter is exactly on such a ray pair, then the entry is $\*$.
Calculating the kneading sequence of any point is substantially
simplified by the observation that, in order to know the entire
kneading sequence at a parameter ray pair of some period $n$, it
suffices to know the first $n-1$ entries in the kneading sequence, so
we only have to look at ray pairs of periods up to $n-1$. This leads
to the following algorithm: for any point $c\in\cz$, find
consecutively the parameter ray pairs of lowest periods between the
previously used ray pair and the point $c$. The periods of these ray
pairs will form a strictly increasing sequence of integers and allow
to reconstruct the kneading sequence, encoding it very efficiently.
If we extend this sequence by a single entry $1$ in the beginning, we
obtain the {\em internal address} of $c$. For details, see
\cite{IntAddr}. In the context of real quadratic polynomials, this
internal address is known as the sequence of {\em cutting times} in
the Hofbauer tower.

The figure shows that certain initial segments of kneading sequences
appear several times. This can be described and explained precisely
and gives rise to certain symmetries of the Mandelbrot set; see
\cite{IntAddr}.

\begin{lemma}[Different Kneading Sequences]
\label{LemKneadingDifferent} \lineclear
Let $c$ be a parabolic parameter and let $z_1$ be the characteristic
periodic point on the parabolic orbit. Among the dynamic rays
landing at $z_1$, only the two characteristic rays can have angles
with identical kneading sequences.
\end{lemma}
\proof
Let $\theta_1,\theta_2,\ldots,\theta_s$ be the angles of the dynamic
rays landing at $z_1$. If their number $s$ is $2$, then both angles
are characteristic, and there is nothing to show. We may hence assume
$s\geq 3$. All the rays $\theta_i$ are periodic of period $n$, say.
By Lemma~\ref{LemRayPermutation}, the orbit period of the parabolic
orbit is exactly $n/s=:k$. Let $z_0$ and $z_0'$ be the two (different)
immediate inverse images of $z_1$ such that $z_0$ is periodic. If any
one of the rays $R(\theta_i)$ is chosen, its two inverse images,
together with any simple path in the critical Fatou component
connecting $z_0$ and $z_0'$, form a partition of the complex plane
into two parts. We label these parts again by $\0$ and $\1$ so that
the dynamic ray at angle $0$ is in part $\0$, and we label the
boundary by $\*$. Now the labels of the parts containing the rays
$R(\theta_i)$, $R(2\theta_i)$, $R(4\theta_i)$, $\ldots$ again
reflect the kneading sequence of $\theta_i$ because the partitions are
bounded at the same angles.

\begin{figure}[tb]
\centerline{\psfig{figure=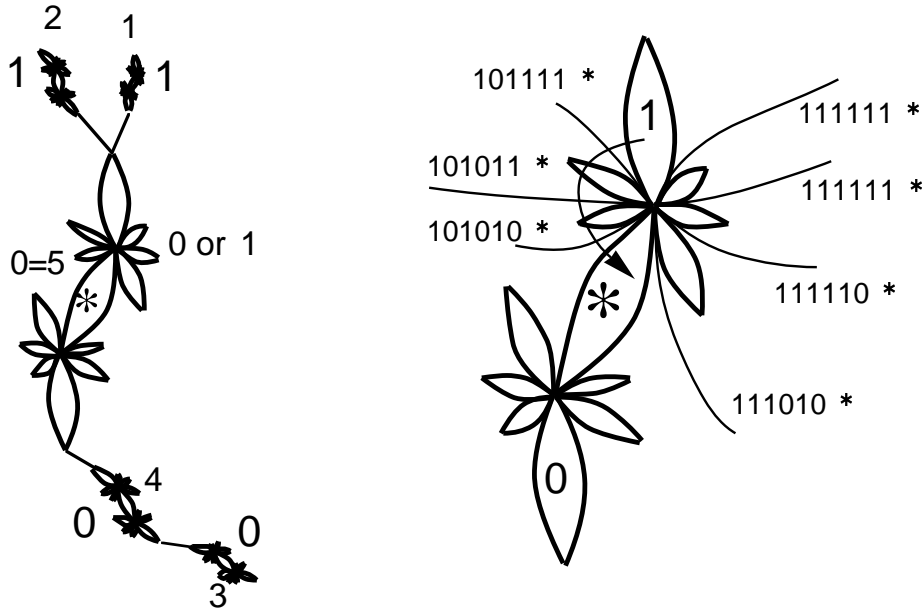,height=82mm}}
\centerline{\parbox{\captionwidth}{ \vspace{0mm} 
\caption[Illustration of the proof of
Lemma~\protect\ref{LemKneadingDifferent}] {\sl Illustration of the
proof of Lemma~\protect\ref{LemKneadingDifferent}. Left: coarse sketch
of the entire Julia set; solid numbers describe the parabolic orbit
(in this case, of period $5$), and outlined numbers specify the
corresponding entries in the kneading sequences of external angles of
rays landing at the characteristic periodic point. Right: blow-up near
the critical point (center, marked by $*$) with the periodic point
$z_0$, considered as a fixed point of the first return map. In this
case, seven rays land at $z_0$ with combinatorial rotation number
$3/7$. The rays are labeled by the corresponding kneading sequences.
The symbols $\0$ and $\1$ indicate regions of the Julia set which are
always on the same side of the partition, independently of which ray
is chosen.}
\label{FigKnSeqDiff}
}}
\end{figure}

Above, we have constructed a tree connecting the parabolic orbit. 
By Lemma~\ref{LemHubbardBranch}, branch points of this tree are on
repelling orbits, so $z_0$ and $z'_0$ have at most two branches of the
tree. One branch always goes into the critical Fatou component to the
critical point. For $z_0$ or $z'_0$, the other branch goes to the
critical value which is always in the region labeled $\1$. The second
branch at the other point ($z'_0$ or $z_0$) must leave in the symmetric
direction, so it will always lead into the region labeled $\0$ (if the
second branch at $z'_0$ were to lead into a direction other than the
symmetric one, then the common image point $z_1$ of $z_0$ and $z'_0$
could not be an endpoint of the tree, which it is by
Lemma~\ref{LemHubbardBranch}). It follows that the entire tree,
except the part between $z_0$ and $z'_0$, is in a subset of the
filled-in Julia set whose label does not depend on which of the angles
$\theta_i$ have been used to define the partition and the kneading
sequence. For positive integers $l$, let $z_l$ be the $l$-th forward
image of $z_0$. If $l$ is not divisible by $k$, it follows that the
label of $z_l$ is independent of $\theta_i$; since
$z_l$ is the landing point of all the dynamic rays at angles
$2^{l-1}\theta_i$, it follows that the $l-1$-st entries of the
kneading sequences of all the $\theta_i$ are the same. Therefore, we
can restrict our attention to the rays at angles
$\theta_1',\theta_2',\ldots,\theta_s'$ landing at $z_0$, where
$\theta_i'$ is an immediate inverse image of $\theta_i$. The first
return dynamics among the angles is multiplication by $2^k$; for the
rays, this must be a cyclic permutation with combinatorial rotation
number $r/s$ for some integer $r$ (Lemma~\ref{LemRayPermutation});
compare Figure~\ref{FigKnSeqDiff}. 

Depending on which of the rays $\theta_i$ is used for the kneading
sequence, i.e. which of the rays $\theta'_i$ defines the partition,
a given ray $\theta'_{j}$ may have label $\0$ or label $\1$. In
particular, the total number among these rays which are in region
$\0$ may be different. But then the number of symbols $\0$ within any
period of the kneading sequence will be different and the kneading
sequences cannot coincide. Two angles among the $\theta_i$ can thus
have the same kneading sequence only if the corresponding partition
has equally many rays in the region labeled $\0$. This leaves only
various pairs of angles at symmetric positions around the critical
value as candidates to have identical kneading sequences. But it is not
hard to verify, looking at the cyclic permutation of the rays
$\theta'_i$, that if two such angles define a partition in which at
least one of the $\theta'_i$ is in region $\0$, then the two
corresponding kneading sequences are different at some position which
is a multiple of $k$. The only two angles with identical kneading
sequences are therefore those for which all the $\theta'_i$ are in
region $\1$ (or on its boundary), so the partition boundary is adjacent
to the Fatou component containing the critical point. The angles
$\theta_i$ are hence those for which the dynamic rays land at
$z_1$ adjacent to the critical value, so the corresponding rays are
the characteristic rays of the parabolic orbit. They do in fact have
identical kneading sequences.
\qed

\proofof{Proposition~\ref{PropAtMostTwoRays}}
We will do the proof by induction on the period $n$. For $n=1$, there
are only two angles $0$ and $1$ which both describe the same parameter
ray, and this ray lands at $c=1/4$. 

To show the statement of the proposition for period $n$, we may suppose
that all parameter rays of periods up to $n-1$ land in pairs. We have a
parabolic parameter $c$ of ray period $n$ and the set $\Theta_c$
contains the angles of parameter rays landing at $c$. All these angles
have the same period $n$ and, by Lemma~\ref{LemKneadingPartition},
identical kneading sequences. Since the corresponding dynamic rays all
land at the characteristic point of the parabolic periodic orbit by
Proposition~\ref{PropNecessary}, Lemma~\ref{LemKneadingDifferent} says
that if $\Theta_c$ contains more than a single element, it contains the
two characteristic angles. This proves the proposition for period $n$,
which implies the Structure Theorem for period $n$, and the inductive
step is completed.
\qed

We have now finished the proof of the periodic part of the theorem,
describing which periodic parameter rays land at common points. This
prepares the ground for combinatorial descriptions such as Lavaurs'
algorithm \cite{La} or internal addresses.

\newsection {Preperiodic Rays}
\label{SecPreperiodic}

In this section, we will turn to parameter rays at preperiodic angles
and show at which Misiurewicz points they land. We will use again
kneading sequences. Recall that if the angle $\theta$ is periodic of
period $n$, then its kneading sequence $\K(\theta)$ will be periodic
of the same period; it will have the symbol $\*$ exactly at positions
$n, 2n, 3n, \ldots$. Moreover, the pointwise limits
$\K_-(\theta):=\lim_{\theta'\nearrow\theta}\K(\theta')$  and
$\K_+(\theta):=\lim_{\theta'\searrow\theta}\K(\theta')$  both exist;
in one of them, all the symbols $\*$ are replaced by $\1$ and in the
other by $\0$ throughout. Both are still periodic; in fact, their
period is $n$ (or a divisor thereof). More precisely, if the parameter
rays at periodic angles $\theta_1$ and $\theta_2$ both land at the
same parameter value, then $\K_\pm(\theta_1)=\K_\mp(\theta_2)$: it
suffices to verify this statement for a single period within the
kneading sequences, and this follows from
Lemma~\ref{LemKneadingPartition}. We can thus imagine every pair of
periodic parameter rays being replaced by two pairs, infinitesimally
close on either side to the given pair and having periodic kneading
sequences without the symbol $\*$.

Kneading sequences of preperiodic rays are themselves preperiodic;
the lengths of the preperiods of external angle and kneading sequence
are equal (this is easy to verify; or see the proof of
Lemma~\ref{LemNumberRaysMisiu}). However, the lengths of the periods
do not have to be equal: the ray $9/56=0.001\,\overline{010}$ has
kneading sequence $\1\1\0\ovl{\1\1\1}=\1\1\0\ovl{\1}$. This fact is
directly related to the number of parameter rays landing at the same
Misiurewicz point; see below.

\begin{figure}[htbp]
\centerline{\psfig{figure=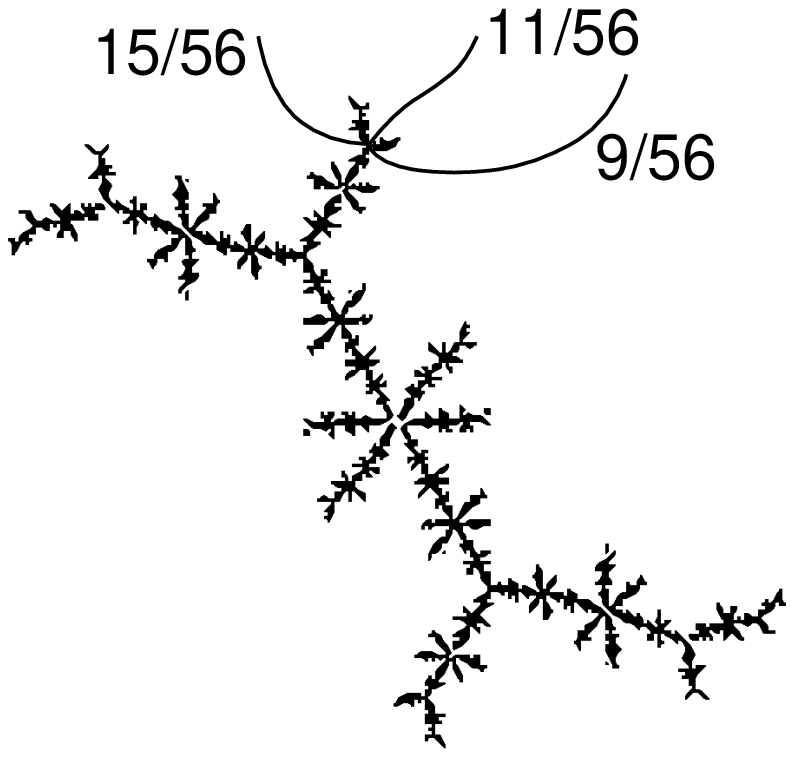,height=50mm} \hfil
	    \psfig{figure=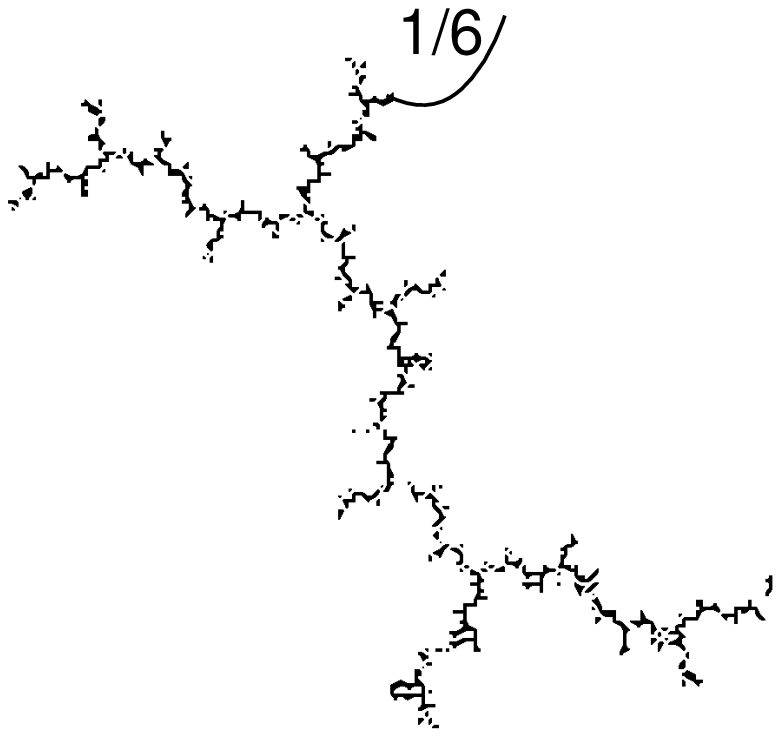,height=50mm}}
\centerline{\parbox{\captionwidth}{ 
\caption{\sl Illustration of the theorem in the preperiodic case. 
Shown are the Julia sets of the polynomials at the landing points of
the parameter rays at angles $9/56$, $11/56$ and $15/56$ (left) and at
angle $1/6$ (right). In both pictures, the dynamic rays landing at the
critical values are drawn. }
\label{FigPreperiodic}
}}
\end{figure}

\vfill

\proofof{Theorem~\ref{ThmRatRays} (preperiodic case)}
Consider any preperiodic parameter ray at angle $\theta$ and let $c$
be one of its limit points. First suppose that $c$ is a parabolic
parameter. We know that there are two parameter rays at periodic angles
$\theta_1,\theta_2$ which land at $c$. We can imagine two parameter ray
pairs infinitesimally close to the ray pair $(\theta_1,\theta_2)$ on
both sides, and these two parameter ray pairs have periodic kneading
sequences without symbols $\*$. Each of these two periodic kneading
sequences must differ at some finite position from the preperiodic
kneading sequence of $\theta$. But we had seen in the previous section
that the regions of constant initial segments of kneading sequences are
bounded by pairs of parameter rays at periodic angles (see
Lemma~\ref{LemKneadingPartition}), so there is a pair of periodic
parameter rays landing at the same point separating the parameter ray
at angle $\theta$ from the two rays at angles $\theta_{1,2}$ and from
the parabolic point $c$. Therefore, $c$ cannot be a limit point of the
parameter ray at angle $\theta$. This contradiction shows that no limit
point of a preperiodic parameter ray is parabolic.

Now we argue similarly as in the proof of
Proposition~\ref{PropPeriodicRaysLand}. For the parameter $c$, there
is no parabolic orbit, so the dynamic ray at angle $\theta$ lands at
a repelling preperiodic point. We want to show that the landing point
is the critical value. Since $c$ is a limit point of the parameter
ray at angle $\theta$, there are arbitrarily close parameters for
which the critical value is on the dynamic ray at angle $\theta$. 

If, for the parameter $c$, the dynamic ray at angle $\theta$ does not
land at the critical point or at a point on the backwards orbit of
the critical point, then the dynamic ray at angle $\theta$ and its
landing point depend continuously on the parameter by
Lemma~\ref{LemStable}, so the critical value must be the landing
point of the dynamic ray at angle $\theta$ for the parameter $c$. If,
however, the landing point of the dynamic $\theta$-ray is on the
backwards orbit of the critical value, then some finite forward image
of this ray will depend continuously on the parameter, and pulling
back may yield a dynamic ray bouncing once into the critical value or
a point on its backward orbit, but after that the two continuations
will land at well-defined points. The ray with both continuations and
both landing points will still depend continuously on the parameter,
so again the dynamic $\theta$-ray must land at the critical value for
the parameter $c$. (However, this contradicts the assumption that the
landing point is on the backwards orbit of the critical value because
that would force the critical value to be periodic.)

We see that, for any limit point $c$ of the parameter ray at angle
$\theta$, the number $c$ is preperiodic under $z\mapsto z^2+c$ with
fixed period and preperiod, and $c$ is a Misiurewicz point. Since any
such point $c$ satisfies a certain polynomial equation, there are only
finitely many such points. The limit set of any ray is connected, so
the parameter ray at angle $\theta$ lands, and the landing point is a
Misiurewicz point with the required properties. This shows the third
part of Theorem~\ref{ThmRatRays}.

For the last part, we have already shown that a Misiurewicz point
cannot be the landing point of a periodic parameter ray, or of a
preperiodic ray with external angle different from the angles of the
dynamic rays landing at the critical value. It remains to show that,
given a Misiurewicz point $c_0$ such that the critical value is the
landing point of the dynamic ray at angle $\theta$, then the
parameter ray at angle $\theta$ lands at $c_0$. We will use ideas
from Douady and Hubbard~\cite{Orsay}. By Lemma~\ref{LemStable}, there
is a simply connected neighborhood $V$ of $c_0$ in parameter space
such that $c_0$ can be continued analytically as a repelling
preperiodic point, yielding an analytic function $z(c)$ with
$z(c_0)=c_0$ such that the dynamic ray at angle $\theta$ for the
parameter $c\in V$ lands at $z(c)$. The relation $z(c)=c$ is
certainly not satisfied identically on all of $V$, so the solutions
are discrete and we may assume that $c_0$ is the only one within $V$. 

Now we consider the winding number of the dynamic ray at angle
$\theta$ around the critical value, which is defined as follows:
denoting the point on the dynamic $\theta$-ray at potential $t\geq 0$
by $z_t$ and decreasing $t$ from $+\infty$ to $0$, the winding number
is the total change of $\arg(z_t-c)$ (divided by $2\pi$ so as to
count in full turns). Provided that the critical value is not on the
dynamic ray or at its landing point, the winding number is
well-defined and finite and depends continuously on the parameter.
When the parameter $c$ moves in a small circle around $c_0$ and if
the winding number is defined all the time, then it must change by an
integer corresponding to the multiplicity of $c$ as a root of
$z(c)-c$. However, when the parameter returns back to where it
started, the winding number must be restored to what it was before.
This requires a discontinuity of the winding number, so there are
parameters arbitrarily close to $c_0$ for which the critical value is
on the dynamic ray at angle $\theta$, and $c_0$ is a limit point of
the parameter ray at angle $\theta$. Since this parameter ray lands,
it lands at $c_0$. This finishes the proof of Theorem~\ref{ThmRatRays}.
\qed

\remark
There is no partition in the dynamic plane showing that preperiodic
parameter rays can not land at parabolic parameters: there are
countably many preperiodic dynamic rays landing at the boundary of
the characteristic Fatou component, for example at preperiodic points
on the parabolic orbit, and they cannot be separated by a stable
partition. 

For the final part of the theorem, we used that a repelling
preperiodic point $z(c)$ depends analytically on the parameter. As
mentioned before, this proof started with a need to describe
parameter spaces of antiholomorphic polynomials like the Tricorn and
Multicorns, and there we do not have analytic dependence on
parameters. Here is another way to prove that every Misiurewicz point
is the landing point of all the parameter rays whose angles are the
external angles of the critical value in the dynamic plane. We start
with any Misiurewicz point $c_0$ and external angle $\theta$ of its
critical value. Let $c_1$ be the landing point of the parameter ray
at angle $\theta$. Then both parameters $c_0$ and $c_1$ have the
property that in the dynamic plane, the ray at angle $\theta$ lands
at the critical value. It suffices to prove that this property
determines the parameter uniquely. This is exactly the content of the
Spider Theorem, which is an iterative procedure to find
postcritically finite polynomials with assigned external angles of
the critical value. In Hubbard and Schleicher~\cite{HS}, there is an
easy proof for polynomials with a single critical point. While the
existence part of that proof works only if the critical point is
periodic, all we need here is the uniqueness part, and that works in
the preperiodic case just as well, both for holomorphic and for
antiholomorphic polynomials.

The last part could probably also be done in a more combinatorial but
rather tedious way, using counting arguments like in the periodic
case. This would, however, be quite delicate, as the number of
Misiurewicz points and the number of parameter rays landing at them
require more bookkeeping: the number of parameter rays landing at
Misiurewicz points varies and can be any positive integer. The
following lemma makes this more precise.

\begin{lemma}[Number of Rays at Misiurewicz Points]
\label{LemNumberRaysMisiu} \lineclear
Suppose that a preperiodic angle $\theta$ has preperiod $l$ and
period $n$. Then the kneading sequence $\K(\theta)$ has the same
preperiod $l$, and its period $k$ divides $n$. If $n/k>1$, then the
total number of parameter rays at preperiodic angles landing at the
same point as the ray at angle $\theta$ is $n/k$; if $n/k=1$, then
the number of parameter rays is $1$ or $2$.
\end{lemma}
In the example above, we had seen that the angle $9/56$ has
period $3$, while its kneading sequence has period $1$. Therefore,
the total number of rays landing at the corresponding Misiurewicz
point is three: their external angles are $9/56,11/56$ and $15/56$.
If more than one ray lands at a given Misiurewicz point, it is not
hard to determine all the angles knowing one of them, using ideas from
the proof below.

\proof
In the dynamic plane of $\theta$, the dynamic ray at angle $\theta$
lands at the critical value, so the two inverse image rays at angles
$\theta/2$ and $(\theta+1)/2$ land at the critical point and separate
the dynamic plane into two parts; this partition cuts the external
angles of dynamic rays in the same way as in the partition defining
the kneading sequence, see Definition~\ref{DefKneadSeq} and
Figure~\ref{FigKneading}. We label the two parts by $\0$ and $\1$ in
the analogous way, assigning the symbol $\*$ to the boundary. The
partition boundary intersects the Julia set only at the critical
point.

The critical value jumps after exactly $l$ steps onto a periodic orbit
of ray period $n$. Denote the critical orbit by $c_0,c_1,c_2,\ldots$
with $c_0=0$ and $c_1=c$, so that $c_{l+1}=c_{l+n+1}$, while
$c_l=-c_{l+n}$. The points $c_l$ and $c_{l+n}$ are on different sides
of the partition. The periodic part of the kneading sequence starts
exactly where the periodic part of the external angles start, so the
preperiods are equal. 

We know that $n$ is the ray period of the orbit the critical value
falls onto. The orbit period is exactly $k$: periodic rays which have
their entire forward orbits on equal sides of the partition land at
the same point, for the following reason: we can connect the landing
points of two such rays by a curve which avoids all the preperiodic
rays landing at the critical point, all the finitely many rays on
their forward orbits, and their landing points (if we have to cross
some of these ray pairs, they must also visit the same sides of the
partition, and we can reduce the problem). Now inverse images of the
rays are connected by inverse images of the curve, which avoids the
same rays. Continuing to take inverse images in this way, the
periodic landing points must converge to each other, so they cannot
be different. 

The number of dynamic rays on the orbit of $\theta$ landing at every
point of the periodic orbit is therefore $n/k$, and the critical value
jumps onto this orbit as a local homeomorphism, so it is the landing
point of equally many preperiodic rays. But these rays reappear in
parameter space as the rays landing at the Misiurewicz point.

It remains to show that there are no extra rays at the periodic
orbit. By Lemma~\ref{LemRayPermutation}, more than two dynamic rays
can land at the same periodic point only if these rays are on the
same orbit, i.e., the dynamics permutes the rays transitively. The
number of dynamic rays can therefore be greater than $n/k$ only
if $n/k=1$, and in that case, there can be at most two rays.
\qed

\remark
It does indeed happen that $n/k=1$ while the number of rays is
two. An example is given by the two parameter rays at angles
$25/56=0/011\,\ovl{100}$ and $31/56=0.100\,\ovl{011}$; their common
kneading sequence is $\1\0\0\,\ovl{\1\0\1}$, so $n=k=3$, but these
two rays land together at a point on the real axis. On the other
hand, for the angle $1/2=0.0\ovl{1}=0.1\ovl{0}$, the kneading
sequence is $\1\ovl{\0}$, so $n=k=1$; the parameter ray at angle
$1/2$ is the only ray landing at the leftmost antenna tip $c=-2$ of
the Mandelbrot set. These rays are indicated in
Figure~\ref{FigMandelRays}.

For a related discussion of rays landing at common points, from the
point of view of ``Thurston obstructions'', see \cite{HS}.

\hide{
At this point, the Farey-algorithm for parameter rays, for kneading
sequences and for the (sub-)limbs which ``hide'' certain cyclic
permutations of dynamic rays might be introduced and used.
}

\newsection{Hyperbolic Components}
\label{SecHyperbolic}

The Orbit Separation Lemma~\ref{LemOrbitSeparation} has an important
consequence: it helps to control the dynamics when a parabolic Julia
set is perturbed. Perturbations of parabolics are a subtle issue
because both the Julia set and the filled-in Julia set behave
drastically discontinuously. We show that nonetheless the landing
points of all the dynamic rays at rational angles behave
continuously wherever the rays land. In a way, the following
proposition is the parabolic analogue to Lemma~\ref{LemStable}, which
dealt with repelling periodic points. However, we will explain below
that the rays themselves do not depend continuously on the parameter. 

\begin{proposition}[Continuous Dependence of Landing Points]
\label{PropContinuousLanding} \lineclear
For any rational angle $\theta$, the landing point of the dynamic ray
at angle $\theta$ depends continuously on the parameter on the entire
subset of parameter space for which the ray lands.
\end{proposition}
The dynamic ray at angle $\theta$ for the polynomial $p_c$ fails to
land if and only if it bounces into the critical point or into a
point on the inverse orbit of the critical point, which happens if
and only if the parameter $c$ is outside the Mandelbrot set on a
parameter ray at one of the finitely many angles
$\{2\theta,4\theta,8\theta,\ldots\}$. All these parameter rays land,
and at the landing parameters, the dynamic ray at angle $\theta$
lands as well. These landing points are the interesting cases of the
proposition.

\proof
First we discuss the case of a periodic angle $\theta$. If the
landing point is repelling, then the proposition reduces to
Lemma~\ref{LemStable}. We may thus assume the landing point to be
parabolic. Under perturbation, any parabolic periodic point splits up
into several periodic points which may be attracting, repelling, or
indifferent, and these periodic points depend continuously on the
parameter. We need to show that the landing point of the ray after
perturbation is one of the continuations of the parabolic periodic
point it was landing at.

Denote the parabolic parameter before perturbation by $c_0$, let $n$
be its ray period and let $V\subset\cz$ be a simply connected open
neighborhood which does not contain further parabolics of equal or
lower ray periods. Then analytic continuation of periodic points of
ray periods up to $n$ is possible in $V-\{c_0\}$. Let $z$ be a
repelling periodic point for the parameter $c_0$; it can then be
continued analytically as a function $z(c)$ in a neighborhood of
$c_0$ so that its orbit remains repelling. We may assume this to be
the case in all of $V$, possibly by shrinking $V$. Then by
Lemma~\ref{LemStable}, $z(c)$ will keep all the rays landing at it
throughout all of $V$; since this point cannot lose rays in $V$, it
cannot gain rays, either. 

Since for the parameter $c_0$, the dynamic ray at angle $\theta$
lands at the parabolic orbit, the landing point of this ray after
perturbation will be on a periodic orbit coming out of the parabolic
orbit, and it remains to show that the landing point does not jump
between continuations of different periodic points of the parabolic
orbit. This is where the Orbit Separation
Lemma~\ref{LemOrbitSeparation} comes in: the parabolic periodic points
are separated by pairs of rays landing at repelling periodic or
preperiodic points, and this separation is stable under perturbations.
The dynamic rays cannot cross this partition, so their landing points
depend continuously on the angle, provided the rays land at all.

For preperiodic rays, the statement follows by taking inverse images
because the pull-back is continuous. If the orbit visits the critical
point along its preperiodic orbit, which happens at Misiurewicz
points, then several preperiodic points may merge and split up with
different rays, but this happens in a continuous way.
\qed
\remark
Unlike their landing points, the dynamic rays themselves may
depend discontinuously on the perturbation. The simplest possible
example occurs near the parabolic parameter $c_0=1/4$: for this
parameter, the dynamic ray at angle $0=1$ lands at the parabolic
fixed point $z=1/2$, and the ray is the real line to the right of
$1/2$. The critical point $0$ is in the interior of the filled-in
Julia set. Perturbing the parameter to the right on the real axis,
i.e., on the parameter ray at angle $0=1$, the dynamic ray will
bounce into the critical point and thus fail to land. But for
arbitrarily small perturbations near this parameter ray, the dynamic
ray at angle $0=1$ will get very close to the critical point before
it turns back and lands near $1/2$. The closer the parameter is to
$c_0=1/4$, the lower will the potential of the critical point be, and
while the dynamic ray keeps reaching out near the critical point, it
does so at lower and lower potentials, and in the limit the part of
the ray at real parts less than $1/2$ will be squeezed off. Points
at any potential $t>0$ will depend continuously on the parameter, and
so does the landing point at potential $t=0$; however, this continuity
is not uniform in $t$, and the dynamic ray as a whole can and does
change discontinuously with respect to the Hausdorff metric. 

Continuous dependence of landing points of rays requires a single
critical point (of possibly higher multiplicity). It is false
already for cubic polynomials; for an example, see the appendix in
Goldberg and Milnor~\cite{GM}.

Most, if not all, of the interior of the Mandelbrot set consists of
what is known as hyperbolic components.
Proposition~\ref{PropContinuousLanding} is one possible key to
understanding many of their properties. First we discuss some
necessary background.

A hyperbolic rational map is one where all the critical points are
attracted by attracting or superattracting periodic orbits. The
dynamical significance is that this is equivalent to the existence of
an expanding metric in a neighborhood of the Julia set, which has
many important consequences such as local connectivity of the Julia
set (see Milnor~\cite{MiIntro}). For a polynomial, the critical point
at $\infty$ is always superattracting, and in the quadratic case, the
polynomial is hyperbolic if the unique finite critical point either
converges to $\infty$ or to a finite (super)attracting orbit.
Hyperbolicity is obviously an open condition. A {\em hyperbolic
component of the Mandelbrot set\/} is a connected component of the
hyperbolic interior. The period of the attracting orbit is constant
throughout the component and defines the {\em period of the
hyperbolic component\/}. We will see below that every boundary point
of a hyperbolic component is a boundary point of the Mandelbrot set,
so a hyperbolic component is also a connected component of the
interior of $\M$. There is no example known of a non-hyperbolic
component; it is conjectured that there are none. A {\em center\/} of
a hyperbolic component is a polynomial for which there is a
superattracting orbit; a {\em root\/} of such a component of period
$n$ is a parabolic boundary point where the parabolic orbit has ray
period $n$. We will show below that every hyperbolic component has a
unique center and a unique root. It is easy to verify that the
multiplier of the attracting orbit on a hyperbolic component is a
proper map from the component to the open unit disk, so it has a
finite mapping degree; we will see that this map is in fact a
conformal isomorphism. The relation between centers and roots of
hyperbolic components is important; the difficulty in establishing it
lies in the discontinuity of Julia sets at parabolic parameters.
Proposition~\ref{PropContinuousLanding} helps to overcome this
difficulty.

\begin{lemma}[Roots of Hyperbolic Components]
\label{LemRoots} \lineclear
Every parabolic parameter with ray period $n$ is a root of at least
one hyperbolic component of period $n$. If the orbit period $k$ is
smaller than $n$, then this parameter is also on the boundary of a
hyperbolic component of period $k$. In no case is such a parabolic
parameter on the boundary of a hyperbolic component of different
period.
\end{lemma}
\proof
First suppose that orbit period and ray period are equal. Then the
first return map of any parabolic periodic point $z$ leaves all the
dynamic rays landing at $z$ fixed, so its multiplier is $+1$. In local
coordinates, the map has the form $\zeta\mapsto
\zeta+\zeta^{q+1}+\ldots$ for some integer $q\geq 1$. The 
point $z$ then has $q$ attracting and repelling petals each, and
every attracting petal must absorb a critical orbit. Since there is a
unique critical point, we have $q=1$. Under perturbation, the
parabolic orbit then breaks up into exactly two orbits of exact
period $n$, and no further orbit is involved. Denote the parabolic
parameter by $c_0$ and let $V$ be a simply connected neighborhood of
$c_0$ not containing further parabolics of equal ray period. In
$V-\{c_0\}$, all periodic points of exact period $n$ can be continued
analytically because their multipliers are different from $+1$. Among
these periodic points, those which are repelling at $c_0$ can be
continued analytically throughout all of $V$, while the two colliding
orbits might be interchanged by a simple loop in $V-\{c_0\}$ (in fact,
they will be: see Corollary~\ref{CorHypBdy}). Their
multipliers are therefore defined on a two-sheeted covering of
$V-\{c_0\}$ and are analytic, even when the point $c_0$ is put back
in. By the open mapping principle, the parameter $c_0$ is on the
boundary of at least one hyperbolic component of period $n$. Since,
for the parameter $c_0$, all the orbits of periods not divisible by
$n$ are repelling, the parameter can be only on the boundary of
hyperbolic components with periods divisible by $n$. If it was on the
boundary of a hyperbolic component with period $rn$ for some integer
$r>1$, then the $rn$-periodic orbit would have to be indifferent at
$c_0$; since there can be only one indifferent orbit, it would have to
merge with the indifferent orbit of period $n$, and this orbit would
get higher multiplicity than $2$, a contradiction. This contradiction
shows that $c_0$ is not on the boundary of any hyperbolic component of
period other than $n$.

If the orbit period strictly divides the ray period, so that
$s:=n/k\geq 2$, then the first return map of the orbit must permute
the rays transitively by Lemma~\ref{LemRayPermutation}. The least
iterate which fixes the rays must also be the least iterate for which
the multiplier is $+1$: the landing point of a periodic ray is either
repelling or has multiplier $+1$ (this is the Snail Lemma, see
\cite{MiIntro}); conversely, whenever the multiplier is $+1$, then all
the finitely many rays must be fixed. It follows that the multiplier of
the first return map of any of the parabolic periodic points is an
exact $s$-th root of unity. The periodic orbit can then be continued
analytically in a neighborhood of the root. Since the multiplier map is
analytic, the parabolic parameter is on the boundary of a hyperbolic
component of period $k$. The $s$-th iterate of the first return map has
multiplier $+1$ and hence again the form $\zeta\mapsto
\zeta+\zeta^{q+1}+\ldots$ in local coordinates, for an integer $q\geq
1$. The number of coalescing fixed points of this iterate is then
exactly $q+1$. 

Since there is only one critical orbit, the first return map of the
parabolic orbit must permute the $q$ attracting petals transitively
and we have $q=s$. For the first, second, \ldots, $s-1$-st iterate of
the first return map, the multiplier is different from $+1$, so the
respective iterate has a single fixed point. The $s$-th iterate,
however, corresponding to the $sk=n$-th iterate of the original
polynomial, has a fixed point of multiplicity $q+1=s+1$: exactly one of
these points has exact period
$k$; all the other points can have no lower periods than $n$, so they
are on a single orbit of period $n$ of which $s$ points each are
coalesced. There is no further orbit involved (or some iterate would
have to have a parabolic fixed point of higher multiplicity with more
attracting petals attached, as above). Since there is a single
indifferent orbit of period $n$, its multiplier is well defined and
analytic in a neighborhood of the parabolic parameter, which is hence
on the boundary of a hyperbolic component of period $n$ as well.
\qed

\smallskip
A root of a hyperbolic component is called {\em primitive\/} if its
parabolic orbit has equal orbit and ray periods, so it is the merger
of two orbits of equal period. If orbit and ray periods are
different, then the root is called {\em non-primitive\/} or a {\em
bifurcation point}: at this parameter, an attracting orbit bifurcates
into another attracting orbit of higher period (the terminology {\em
bi-}furcation comes from the dynamics on the real line, where the
ratio of the periods is always two). We will see below that every
hyperbolic component has a unique root. It therefore makes sense to
call a hyperbolic component primitive or non-primitive according to
whether or not its root is primitive.

\bigskip
We can now draw a couple of useful conclusions.

\goodbreak
\begin{corollary}[Stability at Roots of Hyperbolic Components]\nobreak
\label{CorRootStability} \lineclear
For any hyperbolic component, the landing pattern of periodic and
preperiodic dynamic rays is the same for all polynomials from the
component and at any of its roots.
\end{corollary}
\proof
Again, it suffices to discuss periodic rays; the statement about
preperiodic rays follows simply by taking inverse images because for
the considered parameters, all the preperiodic dynamic rays land, and
they never land at the critical value.
Throughout the component, all the periodic rays land at repelling
periodic points, so no orbit can lose a ray under perturbations, and
consequently no orbit can gain a ray, either. The same rays at
rational angles land at common points throughout the component.

When the dynamics at a primitive root is perturbed into the
component, then the parabolic orbit breaks up into one attracting and
one repelling orbit of equal period, and by continuous dependence of
the landing points (Proposition~\ref{PropContinuousLanding}), the
repelling orbit must inherit all the dynamic rays. It cannot get any
further rays because they have been attached to repelling orbits.

For a non-primitive root, denote orbit and ray periods by $k$ and
$n$, respectively. Perturbing the root into the component of period
$n$, the orbit of period $n$ becomes attracting, and the $k$-periodic
orbit must take all the dynamic rays from the parabolic orbit without
getting any further rays.
\qed
\remark
Perturbing a parabolic orbit with orbit period $k$ and ray period
$n>k$ into the component of period $k$ changes the landing pattern of
rational rays: the parabolic orbit creates an attracting orbit of
period $k$, so a repelling orbit of period $n$ remains, and all the
dynamic rays from the parabolic orbit land at different points.
Phrased differently, when moving from the component of period $k$ into
the one of period $n$, then $n/k$ periodic rays each start landing at
common periodic points; of course, this forces the obvious relations
for the preperiodic rays. The landing pattern of all other periodic
rays remains stable.

This discussion not only describes the landing patterns of periodic
rays within hyperbolic components, but on the entire parameter space
except at parameter rays at periodic angles: since the landing points
depend continuously on the parameter and periodic orbits are simple
except at parabolics, the pattern can change only at parabolic
parameters or at parameter rays where the dynamic rays fail to land. 
The relation between landing patterns of periodic rays and the
structure of parameter space has been investigated and described by 
Milnor in \cite{MiOrbits}. The landing pattern of preperiodic rays
also changes at Misiurewicz points and at preperiodic parameter rays.

\begin{corollary}[The Multiplier Map] 
\label{CorMultiplier} \lineclear
The multiplier map on any hyperbolic component is a conformal
isomorphism onto the open unit disk, and it extends as a
homeomorphism to the closures. In particular, every hyperbolic
component has a unique root and a unique center. The boundary of a
hyperbolic component is contained in the boundary of the Mandelbrot
set.
\end{corollary}
\proof
The multiplier map is a proper analytic map from the hyperbolic
component to the open unit disk and extends continuously to the
boundary. Any point on the boundary has a unique indifferent orbit.
If the multiplier at such a boundary point is different from $+1$,
then orbit and multiplier extend analytically in a neighborhood of
the boundary point. The multiplier is obviously not constant. This
shows in particular that parabolic parameters are dense on the
boundary of the component; since parabolics are landing points of
parameter rays and thus in the boundary of $\M$, the boundary of every
hyperbolic component is contained in the boundary of the Mandelbrot
set. 

The number of parabolic parameters with fixed orbit period and
multiplier $+1$ is finite, so the boundary of any hyperbolic
component consists of a finite number of analytic arcs (which might
contain critical points) limiting on finitely many parabolic
parameters with multipliers $+1$. Since the multiplier map is proper
onto $\disk$, the component has at least one root.

By Corollary~\ref{CorRootStability}, the landing pattern of periodic
dynamic rays has to be the same at all the roots of a given
hyperbolic component. It follows that at all the roots, the angles of
the characteristic rays of the parabolic orbits have to coincide: in
every case, the characteristic ray pair lands at the Fatou component
containing the critical value and separates the critical value from
the rest of the parabolic orbit. Among all ray pairs separating the
critical point from the critical value, the characteristic ray pair
must be the one closest to the critical value. Hence the landing
pattern of periodic dynamic rays determines the characteristic angles.
Since any root of a hyperbolic component must be the landing point of
the parameter rays at the characteristic angles of the parabolic
orbit, every hyperbolic component has a unique root. 

We can now determine the mapping degree $d$, say, of the multiplier
map. The hyperbolic component is simply connected because the
Mandelbrot set is full. Then the multiplier map $\mu$ has exactly
$d-1$ critical points, counting multiplicities. If $d>1$, let
$v\in\disk$ be a critical value of $\mu$ and connect $v$ and $+1$ by
a simple smooth curve $\gamma\subset\disk$ avoiding further critical
values. Then $\mu^{-1}(\gamma)$, together with the unique root of the
component, contains a simple closed curve $\Gamma$ enclosing an open
subset of the hyperbolic component. This subset must map at least
onto all of $\disk-\gamma$, so $\Gamma$ surrounds boundary points of 
the component and thus of $\M$. But this contradicts the fact that
the Mandelbrot set is full, so the multiplier map is a conformal
isomorphism onto $\disk$ and the component has a unique center. It
extends continuously to the closure and is surjective onto
$\partial\disk$ because it is surjective on the component. Near every
non-root of the component, the boundary of the component is an analytic
arc (possibly with critical points) and the multiplier is not locally
constant on the arc, so it is locally injective. Global injectivity now
follows from injectivity on the component. The multiplier map is thus
invertible, and continuity of the inverse is a generality.
\qed

\begin{proposition}[No Shared Roots] 
\label{PropNoSharedRoots} \lineclear
Every parabolic parameter is the root of a single hyperbolic component. 
\end{proposition}
\proof
Since the period of a component equals the ray period of its root, we
can restrict attention to any fixed period $n$. We have well defined
maps from centers to hyperbolic components (which we have just seen
is a bijection) and from hyperbolic components to their roots. This
gives a surjective map from centers of period $n$ to parabolic
parameters of ray period $n$. Denote the number of centers of period
$n$ by $s_n$.

A center of a hyperbolic component of period $n$ is a point $c$ such
that $0$ is periodic of exact period $n$ under $z\mapsto z^2+c$;
therefore, $c$ must satisfy a polynomial equation
$(\ldots((c^2+c)^2+c)\ldots)^2+c=0$ of degree $2^{n-1}$. Since this
polynomial is also solved by centers of components of periods $k$
dividing $n$, we get the recursive relation $\sum_{k|n} s_k=2^{n-1}$.
By Lemma~\ref{LemParaCount}, this is exactly the number of parabolic
parameters of ray period $k$. Since a surjective map between finite
sets of equal cardinality is a bijection, every parabolic parameter is
the root of a single hyperbolic component.
\qed
\remark
This proposition shows even without resorting to
Corollary~\ref{CorMultiplier} that every hyperbolic component has a
unique center, so that the only critical point of the multiplier map
(if it had mapping degree greater than one) could be the center. This
is indeed what happens for the ``Multibrot sets'': the connectedness
loci for the maps $z\mapsto z^d+c$ with $d\geq 2$.

Before continuing the study of hyperbolic components, we note an
algebraic observation following from the proof we have just given.
\begin{corollary}[Centers of Components as Algebraic Numbers]
\label{CorCentersAlgebraic} \lineclear 
Every center of a hyperbolic component of degree $n$ is an algebraic
integer of degree at most $s_n$. It is a simple root of its minimal
polynomial.
\qedd
\end{corollary}

A neat algebraic proof for this fact has been given by Gleason; see
\cite{Orsay}. As far as I know, the algebraic structure of the
minimal polynomials of the centers of hyperbolic components is not
known: when factored according to exact periods, are they
irreducible? What are their Galois groups? Manning (unpublished) has
verified irreducibility for $n\leq 10$, and he has determined that
the Galois groups for the first few periods are the full symmetric
groups. Giarrusso (unpublished) has observed that this induces a
Galois action between the Riemann maps of hyperbolic components of
equal periods, provided that their centers are algebraically
conjugate.

Now we can describe the boundary hyperbolic components much more
completely.
\begin{corollary}[Boundary of Hyperbolic Components]
\label{CorHypBdy} \lineclear
No non-parabolic parameter can be on the boundary of more than one
hyperbolic component. Every parabolic parameter is either a primitive
root of a hyperbolic component and on the boundary of no further
component, or it is a ``point of bifurcation'': a non-primitive root
of a hyperbolic component and on the boundary of a unique further
hyperbolic component. In particular, if two hyperbolic components
have a boundary point in common, then this point is the root of
exactly one of them.  The boundary of a hyperbolic component is a
smooth analytic curve, except at the root of a primitive component.
At a primitive root, the component has a cusp, and analytic
continuation of periodic points along a small loop around this cusp
interchanges the two orbits which merge at this cusp.
\end{corollary}
\proof
If two hyperbolic components have a non-parabolic parameter in their
common boundary, then the landing patterns of periodic rays must be
the same within both components. This must then also be true at their
respective roots, which yields a contradiction: on the one hand, the
roots must be different by Proposition~\ref{PropNoSharedRoots}; on
the other hand, the parabolic orbits at the roots must have the same
characteristic angles (compare the proof of
Corollary~\ref{CorMultiplier}), so they must be the landing points of
the same parameter rays. It follows that the multiplier map of the
indifferent orbit cannot have a critical point at $c_0$: if it had a
critical point there, then $c_0$ would connect locally two regions of
hyperbolic parameters which cannot belong to different hyperbolic
components; however, if they belonged to the same component, then the
closure of the component would separate part of its boundary from the
exterior of the Mandelbrot set, a contradiction. Therefore, the
boundary of every hyperbolic component is a smooth analytic curve
near every non-parabolic boundary point.

Now let $c_0$ be a parabolic parameter of ray period $n$ and orbit
period $k$. We know that it is the root of a unique hyperbolic
component of period $n$. 

In the non-primitive case (when $k$ strictly divides $n$), the point
$c_0$ cannot be on the boundary of a hyperbolic component of period
different from $n$ and $k$ by Lemma~\ref{LemRoots}. It is on the
boundary of a single hyperbolic component of period $n$
(Proposition~\ref{PropNoSharedRoots}). In a small punctured
neighborhood of $c_0$ avoiding further parabolics of ray period $n$,
the multiplier map of the $n$-periodic orbit is analytic and cannot
have a critical point at $c_0$, for the same reason as above, so the
component of period $n$ occupies asymptotically (on small scales) a
half plane near $c_0$. The parameter $c_0$ is also on the boundary
of a component of period $k$ the multiplier of which is analytic near
$c_0$. Since hyperbolic components cannot overlap, this component
must asymptotically be contained in a half plane, so the multiplier
map cannot have a critical point and must then be locally injective
near $c_0$. The boundaries of both components must then be smooth
analytic curves near $c_0$.

In the primitive case $k=n$, the parameter $c_0$ cannot be on the
boundary of a hyperbolic component of period different from $n$ by
Lemma~\ref{LemRoots}, and it cannot be on the boundary of two
hyperbolic components of period $n$ because otherwise it would have to
be their simultaneous root, contradicting
Proposition~\ref{PropNoSharedRoots}. In a small simply connected
neighborhood $V$ of $c_0$, analytic continuation of the two orbits
colliding at $c_0$ is possible in $V-\{c_0\}$ (compare the proof of
Lemma~\ref{LemRoots}), so their multipliers can be defined on a
two-sheeted cover of $V$ ramified at $c_0$. If analytic continuation
of these two orbits along a simple loop in $V$ around $c_0$ did not
interchange the two orbits, then both multipliers could be defined in
$V$, and both would define different hyperbolic components
intersecting $V$ in disjoint regions, yielding the same contradiction
as in the non-primitive case above. Therefore, small simple loops
around $c_0$ do interchange the two orbits. In order to avoid the
same contradiction again, the multiplier must be locally injective on
the two-sheeted covering on $V$. Projecting down onto $V$, the
component must asymptotically occupy a full set of directions, so the
component has a cusp.
\qed

\remark
The basic motor for many of these proofs about hyperbolic components
was uniqueness of parabolic parameters with given combinatorics, via
landing properties of parameter rays at periodic angles. A
consequence was uniqueness of centers of hyperbolic components with
given combinatorics. One can also turn this discussion around and
start with centers of hyperbolic components: the fact that hyperbolic
components must have different combinatorics, and that they have
unique centers, is a consequence of Thurston's topological
characterization of rational maps, in this case most easily used in
the form of the Spider Theorem \cite{HS}.

\newpage


\bigskip
\noindent

\hfill
\parbox[t]{65mm}{
Dierk Schleicher \\
Fakult\"at f\"ur Mathematik\\
Technische Universit\"at\\
Barer Stra{\ss}e 23\\
D-80290 M\"unchen, Germany\\
{\sl dierk$@$mathematik.tu-muenchen.de}
}


\end{document}